\documentclass{article}

\usepackage{arxiv}

\usepackage[utf8]{inputenc} 
\usepackage[T1]{fontenc}    
\usepackage{hyperref}       
\usepackage{url}            
\usepackage{booktabs}       
\usepackage{amsfonts}       
\usepackage{nicefrac}       
\usepackage{microtype}      
\usepackage{lipsum}		
\usepackage{graphicx}
\usepackage{natbib}
\usepackage{doi}
\usepackage{amsmath}

\title{Exploring the origins of switching dynamics in a multifunctional reservoir computer}

\author{ \href{https://orcid.org/0000-0002-4968-8972}{\includegraphics[scale=0.06]{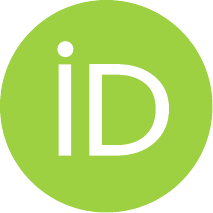}\hspace{1mm}Andrew Flynn}\thanks{Corresponding author.} \\
	School of Mathematical Sciences, University College Cork, Cork T12 XF62, Ireland,\\ INFANT Research Centre, University College Cork, T12 DC4A Cork, Ireland.\\
	\texttt{andrew\_flynn@umail.ucc.ie} \\
	\And
	\href{https://orcid.org/0000-0002-5865-3712}{\includegraphics[scale=0.06]{orcid.pdf}\hspace{1mm}Andreas Amann} \\
	School of Mathematical Sciences, University College Cork, Cork T12 XF62, Ireland,\\
 Potsdam Institute for Climate Impact Research, Telegrafenberg, 14473 Potsdam, Germany.\\
	\texttt{a.amann@ucc.ie} \\
}

\renewcommand{\shorttitle}{Origins of switching dynamics in a multifunctional reservoir computer}

\begin{document}
\maketitle

\begin{abstract}
The concept of multifunctionality has enabled reservoir computers (RCs), a type of dynamical system that is typically realised as an artificial neural network, to reconstruct multiple attractors simultaneously using the same set of trained weights. However there are many additional phenomena that arise when training a RC to reconstruct more than one attractor. Previous studies have found that, in certain cases, if the RC fails to reconstruct a coexistence of attractors then it exhibits a form of metastability whereby, without any external input, the state of the RC switches between different modes of behaviour that resemble properties of the attractors it failed to reconstruct. In this paper we explore the origins of these switching dynamics in a paradigmatic setting via the `seeing double' problem.
\end{abstract}


\section{Introduction}

Multifunctionality is the term used to describe a neural network that has the ability to perform multiple tasks without changing any of its connections.
Multifunctionality is an essential property of certain biological neural networks and has been an active area of research in neuroscience since the mid-1980s with seminal work published in \cite{Mpitsos86MFNN} and in \cite{getting89Principles} followed by further review papers by \cite{Dickinson95MF}, \cite{Marder96MFprinciples}, and more recently reviewed in \cite{briggman08multifunctional}. These studies have identified that a multifunctional neural networks in principle resembles a multistable dynamical system. In this sense, for each task that the network performs there is an attractor associated with it. This attractor is in coexistence with several other attractors in the network's state space and each attractor is related to one of the tasks that the network performs. Therefore, in order to perform a given task, the multifunctional network requires a cue in the form of an initial condition in the basin of attraction of the attractor associated with the task. 

Taking all of the above into account, where this ability to harness multistability becomes of immediate relevance is in the domain of machine learning (ML) as multifunctionality can be used to unlock additional computational capabilities of artificial neural networks (ANNs) that would otherwise have remained dormant. In \cite{flynn2021multifunctionality} multifunctionality was achieved in an artificial setting for the first time via the reservoir computing approach to ML. This involved training a `reservoir computer' (RC), which in this case was a dynamical system in the form of an artificial neural network, to reconstruct a coexistence of chaotic attractors from different dynamical systems using the same set of trained weights. This RC was driven with input from these chaotic attractors and the RC's response dynamics to the different driving inputs were used to obtain a readout layer to replace the drive, after which the RC became a multistable system that reconstructed a coexistence of the chaotic attractors. In this example, to perform a particular task, i.e. to reconstruct a particular chaotic attractor, the multifunctional RC is like any other multistable dynamical system and needs only to be initialised with an initial condition in the basin of attraction of the corresponding attractor. 

There are many additional phenomena that can arise and factors to consider when training a RC to reconstruct more than one attractor simultaneously. For instance, it was shown in \cite{flynn2021multifunctionality} that multifunctionality becomes increasingly difficult to achieve the greater the difference of the time scales of the attractors that the RC is trained to reconstruct.
Furthermore, in \cite{flynn2023seeingdouble} where the RC was trained to solve the 'seeing double' problem which involves training the RC to construct a coexistence of attractors that describe clockwise and anticlockwise trajectories on two circular orbits, it was shown that by manually shifting the location of the training data describing these orbits, the closer the orbits are to one another the more difficult it is for the RC to achieve multifunctionality. Remarkably, for a small range of training parameters, it was found that the RC achieves multifunctionality even when the orbits are overlapping in state space (in the sense that the training data used to drive the RC contains identical data points from the different orbits). 
In \cite{flynn2023seeingdouble} and \cite{flynn2023PhD}, it was shown that in certain cases when the RC fails to achieve multifunctionality it instead produces a variety of episodic switching patterns between different metastable states that resemble the dynamics it failed to reconstruct. From further investigation of the seeing double problem we have found a similar phenomenon to occur when the orbits are moved closer together. The purpose of this paper is to examine the origins of these switching dynamics in much greater detail.


In this paper we explore the origins of the transition from multifunctionality to metastable switching dynamics in much greater detail. We find that for a small change in the spectral radius of the RC's internal connectivity matrix, the RC first fails to reconstruct one of the orbits as the corresponding reconstructed attractor becomes unstable and it is only after a relatively long transient that the RC approaches the other reconstructed orbit (which is the only stable attractor present in the system). After another small change in the spectral radius, the other reconstructed orbit also becomes unstable, this results in the state of the RC switching between the dynamics of these two unstable states. From closer inspection we find that when the second attractor becomes unstable there is a new attractor created that facilitates these switching dynamics. This new attractor is created through this sequence of attractors becoming unstable because due to the RC's design it is prohibited from becoming globally unstable. We show that these switching dynamics appear when the orbits are brought closer together, touch, and overlap. 
From computing the probability density of different residence times in each of the metastable states, we find a sawtooth like pattern consisting of multiple branches of exponentially distributed points where each branch describes a particular path taken by the RC on each of the metastable states.

\section{Methods}

In this section we introduce the particular RC that is studied throughout this paper, we describe how this RC is trained to achieve multifunctionality, and we outline the specifics of the seeing double problem, the task that the RC is trained to solve. We follow the same procedure as in \cite{flynn2023seeingdouble}.

\subsection{\label{sec:RCintro}Reservoir computing}

\subsubsection{Central philosophy of reservoir computing}

Today the term `reservoir computer' (RC) is generally used to describe a dynamical system that, for instance, can be realised as an artificial neural network (ANN) and trained to solve certain machine learning (ML) problems without explicitly training the internal structure of the system. As outlined in \cite{nakajima2021RCbook}, the reservoir computing approach to ML received its name in \cite{verstraeten2005firstRCmention} where the term reservoir computer was coined as a means to establish a new ML framework based on the common concepts of Echo-State Networks \cite{Jaeger01ESN} (ESNs) and Liquid-State Machines \cite{Maass02_LSM} (LSMs). These are two independently proposed designs of ANNs with recurrent connections (RNNs) that share the following philosophy; instead of training all the weights in a network it is sufficient to only optimise the weights of a readout layer in order to solve a particular problem. This ideological shift in training RNNs stems from the design of a suitable internal layer, known as the `reservoir', which does not need to be trained according to a given task. The role of this reservoir is to enable the state of the RC to become a representation of the history of training input signals related to a particular task, only a readout layer needs to be found in order to project this information out of the RC to solve the given task. Multifunctionality extends this philosophy by demonstrating that a RC's response to several different sequences of training input signals, each of which related to a particular task, can be harnessed to produce a single RC that performs all of these tasks using the same readout layer. 


\subsubsection{RC formulation}\label{sssec:OpenLoopRC}

The RC that is studied throughout this paper was introduced in \cite{LuHuntOtt18RC}, before this RC is trained it is defined as the following ANN in the form of a nonautonomous dynamical system which we refer to as the `open-loop RC', 
\begin{align}
    \dot{\boldsymbol{r}}(t) &= \gamma \left[ - \boldsymbol{r}(t) + \tanh{\left( \textbf{M} \, \boldsymbol{r}(t) + \sigma \textbf{W}_{in} \, \boldsymbol{u}(t) \right)} \right], \label{eq:ListenRes}\\
    \boldsymbol{r}\left( 0 \right) &= \boldsymbol{0}^{T}.\label{eq:ListenResIC}
\end{align}
In Eq.\,\eqref{eq:ListenRes}, $\boldsymbol{r}(t) \in \mathbb{R}^{N}$ describes the state of the open-loop RC at a given time $t$ and $N$ is the number of artificial neurons in the network. Solutions of Eq.\,\eqref{eq:ListenRes} are computed using the 4$^{th}$ order Runge-Kutta method with time step $\tau = 0.01$. $\gamma$ is a decay-rate parameter that arises during the derivation of Eq.\,\eqref{eq:ListenRes} from the discrete-time design proposed by \cite{Jaeger01ESN}. The $\tanh$ `activation function' is a pointwise operation and is defined as $\tanh\left( \cdot \right) : \mathbb{R}^{N} \to \mathbb{R}^{N}$. The adjacency matrix, $\textbf{M} \in \mathbb{R}^{N \times N}$, plays the role of the `reservoir'. The input strength parameter, $\sigma$, and the input matrix, $\textbf{W}_{in} \in \mathbb{R}^{N \times D}$, when multiplied together represent the weight given to the $D$-dimensional driving input, $\boldsymbol{u}(t) \in \mathbb{R}^{D}$, as it is projected into the open-loop RC. We use the superscript $T$ to denote the vector transpose operation.

The elements of $\textbf{M}$ and $\textbf{W}_{in}$ are the same that are used in \cite{flynn2023seeingdouble} in order to provide a direct comparison to the results of this present paper. This $\textbf{M}$ was designed by first constructing a random sparse matrix where each element is chosen independently to be nonzero with probability $P=0.04$ (i.e. sparsity $= P$ or degree $= N/P$) and these nonzero elements are chosen uniformly from $\left( -1, 1 \right)$. The elements of this random sparse matrix are subsequently rescaled so that the resultant matrix, that we then call $\textbf{M}$, has a specific spectral radius denoted by $\rho$, which is the magnitude of the largest eigenvalue of $\textbf{M}$. The corresponding input matrix, $\textbf{W}_{in}$, was designed such that each row has only one nonzero randomly assigned element that is chosen uniformly from $\left( -1, 1 \right)$. As a result each neuron is driven with only one component of $\boldsymbol{u}(t)$.

Building on the results of \cite{flynn2023seeingdouble}, $\rho$ is again shown to play a significant role in producing the switching dynamics that are studied in this paper. $\rho$ has also been a key parameter in previous results on training a RC to achieve multifunctionality, see: \cite{flynn2021multifunctionality,flynn2021symmetry,Flynn22_LimitsMF,morraflynn23_MF_fly}. One of the main reasons why $\rho$ is such an influential parameter of this RC is that it is used to tune how previous states of the RC impact the current state. This becomes particularly important in scenarios involving overlapping training data because there must be a sufficiently large weight placed on previous states in order for the RC to distinguish between identical data points from the different sets of training data. In this paper we find that when $\rho$ is not sufficiently large then the RC cannot easily distinguish between the different orbits which, in certain scenarios, leads to the state of the RC switching between the orbits.


\subsubsection{Training a RC to achieve multifunctionality}\label{sssec:RCtraining}

We now outline the steps involved in training Eq.~\eqref{eq:ListenRes} to achieve multifunctionality. To illustrate this procedure, we consider the case of training the open-loop RC in Eq.~\eqref{eq:ListenRes} to reconstruct a coexistence of two attractors, $\mathcal{A}_{1},\, \mathcal{A}_{2}, \subset \mathbb{R}^{D}$, given access to a trajectory on each attractor described by $\boldsymbol{u}_{\left(\mathcal{A}_{1}\right)}(t) \in \mathcal{A}_{1}$ and $\boldsymbol{u}_{\left(\mathcal{A}_{2}\right)}(t) \in \mathcal{A}_{2}$. In the case of multifunctionality, the aim of the training is to determine a `readout function/layer', defined as $\hat{\boldsymbol{\psi}}\left( \cdot \right): \mathbb{R}^{2 N} \to \mathbb{R}^{D}$, that enables us to replace $\boldsymbol{u}(t)$ in Eq.\,\eqref{eq:ListenRes} with $\hat{\boldsymbol{\psi}}\left( \cdot \right)$ and form a new `closed-loop RC' which is capable of reconstructing a coexistence of $\mathcal{A}_{1}$ and $\mathcal{A}_{2}$. In this paper $\hat{\boldsymbol{\psi}}\left( \cdot \right)$ is constructed as, 
\begin{align}
\hat{\boldsymbol{\psi}}\left(\boldsymbol{r}(t)\right) = \textbf{W}_{out} \boldsymbol{q}( \boldsymbol{r}(t) ),\label{eq:ReadoutFunction}
\end{align}
where $\textbf{W}_{out} \in \mathbb{R}^{D \times 2N}$ is the `readout matrix' and $\boldsymbol{q}( \boldsymbol{r}(t) ) \in \mathbb{R}^{2 N}$ is given by,
\begin{align}
    \boldsymbol{q}(\boldsymbol{r}(t))=\left( \boldsymbol{r}(t) \,\,
\boldsymbol{r}^{2}(t) \right)^{T},\label{eq:q_square}
\end{align}
where $\boldsymbol{r}^{2}(t) = \left( r_{1}^{2}(t), r_{2}^{2}(t), \ldots, r_{N}^{2}(t) \right)^{T}$. The purpose of $\boldsymbol{q}(\cdot)$ is to prevent the occurrence of `mirror-attractors' which can impede the ability of the RC to reconstruct attractors as reported on in \cite{herteux2020Symm} and \cite{flynn2021symmetry}. To compute $\textbf{W}_{out}$ in Eq.~\eqref{eq:ReadoutFunction} we use a ridge regression technique, this consists of solving the following equation,
\begin{align}
    \textbf{W}_{out} = \textbf{Y}_{C} \textbf{X}_{C}^{T} \left( \textbf{X}_{C} \textbf{X}_{C}^{T} + \beta \, \textbf{I} \right)^{-1},\label{eq:WoutRegression}
\end{align}
where $\beta$ is the `ridge parameter' and is tuned to reduce the magnitudes of elements in $\textbf{W}_{out}$ in order to discourage overfitting, $\textbf{I}$ is the identity matrix, and $\textbf{X}_{C}$ and $\textbf{Y}_{C}$ are the training data matrices which are both constructed as concatenations of two smaller matrices where $\textbf{X}_{C}=\left[ \textbf{X}_{\mathcal{A}_{1}}, \textbf{X}_{\mathcal{A}_{2}} \right]$ and $\textbf{Y}_{C}=\left[ \textbf{Y}_{\mathcal{A}_{1}}, \textbf{Y}_{\mathcal{A}_{2}} \right]$. The elements of these $\textbf{X}_{\mathcal{A}_{1}}$ and $ \textbf{X}_{\mathcal{A}_{2}}$ matrices are computed as follows: we first drive the open-loop RC in Eq.\,\eqref{eq:ListenRes} with input $\boldsymbol{u}_{\left(\mathcal{A}_{1}\right)}(t) \in \mathcal{A}_{1}$ for $0 < t \leq t_{\text{train}}$ and then repeat this process for $\boldsymbol{u}_{\left(\mathcal{A}_{2}\right)}(t) \in \mathcal{A}_{2}$. The corresponding responses of the open-loop RC to these driving inputs are denoted by $\boldsymbol{r}_{\left(\mathcal{A}_{1}\right)}(t)$ and $\boldsymbol{r}_{\left(\mathcal{A}_{2}\right)}(t)$. It is these responses that are used to generate the elements of $\textbf{X}_{\mathcal{A}_{1}}$ and $ \textbf{X}_{\mathcal{A}_{2}}$ where
\begin{align}
    \textbf{X}_{\mathcal{A}_{1}} = \left[ \begin{array}{cccc}
    \boldsymbol{q}(\boldsymbol{r}_{\left(\mathcal{A}_{1}\right)}(t_{\text{listen}})) & \boldsymbol{q}(\boldsymbol{r}_{\left(\mathcal{A}_{1}\right)}(t_{\text{listen}}+\tau)) 
    &
    \cdots 
    &
    \boldsymbol{q}(\boldsymbol{r}_{\left(\mathcal{A}_{1}\right)}(t_{\text{train}}))
    \end{array} \right],\label{eq:Xmat}
\end{align}
and similarly for $\textbf{X}_{\mathcal{A}_{2}}$. The elements of the corresponding $\textbf{Y}_{\mathcal{A}_{1}}$ and $ \textbf{Y}_{\mathcal{A}_{2}}$ matrices are defined as,
\begin{align}
    \textbf{Y}_{\mathcal{A}_{1}} = \left[ \begin{array}{cccc}
    \boldsymbol{u}_{\left(\mathcal{A}_{1}\right)}(t_{\text{listen}}) & \boldsymbol{u}_{\left(\mathcal{A}_{1}\right)}(t_{\text{listen}}+\tau)
    &
    \cdots 
    &
    \boldsymbol{u}_{\left(\mathcal{A}_{1}\right)}(t_{\text{train}})
    \end{array} \right],\label{eq:Ymat}
\end{align}
and similarly for $\textbf{Y}_{\mathcal{A}_{2}}$. The time $t_{\text{listen}}$ is chosen such that at this time both $\boldsymbol{r}_{(\mathcal{A}_{1})}(t)$ and $\boldsymbol{r}_{(\mathcal{A}_{2})}(t)$ are determined by a history of driving inputs and are no longer dependent on the open-loop RC's initial condition, the duration of time from $t=0$ to $t=t_{\text{listen}}$ is known as `the listening stage'. The time $t_{\text{train}}$ is chosen such that $\textbf{X}_{\mathcal{A}_{1}}$ and $ \textbf{X}_{\mathcal{A}_{2}}$ contain a sufficiently long representation of a trajectory on $\mathcal{A}_{1}$ and $\mathcal{A}_{2}$, the duration of time from $t=t_{\text{listen}}$ to $t=t_{\text{train}}$ is known as `the training stage'. It is important to highlight that $\textbf{M}$, $\textbf{W}_{in}$, and all training parameters remain identical when generating $\textbf{X}_{\mathcal{A}_{1}}$ and $ \textbf{X}_{\mathcal{A}_{2}}$.

\subsubsection{The `closed loop' RC}\label{sssec:ClosedLoopRC}

After following the steps outlined in the previous section and obtaining $\textbf{W}_{out}$ from Eq.~\eqref{eq:WoutRegression}, $\boldsymbol{u}(t)$ in Eq.\,\eqref{eq:ListenRes} can then be replaced by $\hat{\boldsymbol{\psi}} \left( \boldsymbol{r}(t) \right)$. In Eq.\,\eqref{eq:PredRes} we now define the resulting closed-loop RC as the following autonomous dynamical system,
\begin{align}
    \hspace{-0.1cm}\dot{\hat{\boldsymbol{r}}}(t) &= \gamma \left[ - \hat{\boldsymbol{r}}(t) + \tanh{\left( \textbf{M} \, \hat{\boldsymbol{r}}(t) + \sigma \textbf{W}_{in}  \textbf{W}_{out} \, \boldsymbol{
    q}(\hat{\boldsymbol{r}}(t)) \right)} \right], \label{eq:PredRes}\\
    \hspace{-0.1cm}\hat{\boldsymbol{r}}\left( 0 \right) &= \boldsymbol{r}\left( t_{\text{train}} \right),\label{eq:PredResIC}
\end{align}
where $\hat{\boldsymbol{r}}(t)$ denotes the state of the closed-loop RC at a given time $t$. While $\hat{\boldsymbol{r}}(t)$ and $\boldsymbol{r}(t)$ are both $N$-dimensional vectors, the purpose of this notation is to distinguish between the dynamics of the closed-loop and open-loop RCs. Furthermore, we consider $\hat{\boldsymbol{r}}(t) \in \mathbb{S}$ where $\mathbb{S}$ is referred to as the `RC's state space' and is used henceforth when discussing the dynamics of the closed-loop RC in $\mathbb{R}^{N}$. By computing solutions of Eq.\,\eqref{eq:PredRes}, predictions of $\boldsymbol{u}(t)$ for $t>t_{\text{train}}$, denoted as $\hat{\boldsymbol{u}}(t)$, are given by,
\begin{align}
    \hat{\boldsymbol{u}}(t) = \hat{\boldsymbol{\psi}}\left( \hat{\boldsymbol{r}}(t)\right).\label{eq:uhat}
\end{align}
Again, while both $\boldsymbol{u}(t)$ and $\hat{\boldsymbol{u}}(t)$ are $D$-dimensional vectors, we use the same convention to indicate that $\hat{\boldsymbol{u}}(t)$ is a prediction of $\boldsymbol{u}(t)$ at time $t$. We also define $\hat{\boldsymbol{u}}(t) \in \mathbb{P}$ where $\mathbb{P}$ is referred to as the `projected state space' and is used henceforth when discussing these projected dynamics of the closed-loop RC.



To test whether the closed-loop RC has achieved multifunctionality, we initialise Eq.~\eqref{eq:PredRes} with $\hat{\boldsymbol{r}}\left( 0 \right) = \boldsymbol{r}_{\left(\mathcal{A}_{1}\right)}\left( t_{train} \right)$ and $\boldsymbol{r}_{\left(\mathcal{A}_{2}\right)}\left( t_{train} \right)$ and from these initial conditions we examine the long-term projected dynamics of the closed-loop RC in $\mathbb{P}$. We say the closed-loop RC has achieved multifunctionality if the long-term dynamical characteristics of $\hat{\boldsymbol{u}}_{(\mathcal{A}_{1})}(t)$ and $\hat{\boldsymbol{u}}_{(\mathcal{A}_{2})}(t)$ are indistinguishable from $\boldsymbol{u}_{(\mathcal{A}_{1})}(t)$ and $\boldsymbol{u}_{(\mathcal{A}_{2})}(t)$. If this is the case then we can say that there exists a coexistence of attractors $\mathcal{S}_{1}, \mathcal{S}_{2} \subset \mathbb{S}$ and when the state of the closed-loop RC approaches either $\mathcal{S}_{1}$ or $\mathcal{S}_{2}$, the corresponding projected dynamics in $\mathbb{P}$ are referred to as the `reconstructed attractors', $\hat{\mathcal{A}}_{1}, \hat{\mathcal{A}}_{2} \subset \mathbb{P}$, which resemble the long-term dynamics of $\mathcal{A}_{1}$ and $\mathcal{A}_{2}$. By resembling the long-term dynamics it is meant that, for instance, $\mathcal{A}_{1}$ and $\hat{\mathcal{A}}_{1}$ will have nearly identical Poincar{\'e} sections when computed for the same region of $\mathbb{R}^{D}$ and $\mathbb{P}$ as $t \to \infty$. If multifunctionality is achieved then we refer to the resulting multistable closed-loop RC as the `multifunctional RC'. 

We comment that $\hat{\boldsymbol{r}}\left( 0 \right) = \boldsymbol{r}_{\left(\mathcal{A}_{1}\right)}\left( t_{\text{train}} \right)$ and $\boldsymbol{r}_{\left(\mathcal{A}_{2}\right)}\left( t_{\text{train}} \right)$ are not the only initial conditions that will allow the closed-loop RC to reconstruct $\mathcal{A}_{1}$ and $\mathcal{A}_{2}$, so long as the closed-loop RC is initialised with a point in the basin of attraction of either $\mathcal{S}_{1}$ or $\mathcal{S}_{2}$ then the corresponding attractor will be reconstructed in $\mathbb{P}$.

\subsection{Seeing Double}\label{ssec:SD_sec}

The specifics of the `seeing double' problem are outlined in this section. This numerical experiment was introduced in \cite{flynn2023seeingdouble} as a means to systematically study the issues related to multifunctionality and overlapping training data.

\subsubsection[Numerical experiment setup]{Numerical experiment setup}\label{sec:SD_Experiment}

The seeing double problem consists of training a RC to construct a coexistence of attractors such that their dynamics in $\mathbb{P}$ follow trajectories along two circular orbits of equal radius and rotate in opposite directions. The difficulty of this task is varied by moving the centres of these orbits closer together or further apart. When these orbits are overlapping, the RC is therefore required to distinguish between points in $\mathbb{R}^{D}$ that are common to both cycles in order to exhibit multifunctionality.

The driving input to the RC is generated via,
\begin{equation}
\boldsymbol{u}(t)=
    \left( \begin{array}{c}
        x(t) \\
        y(t)
    \end{array} \right)
    = \left( \begin{array}{c}
        b_{x} \cos{\left( t \right)} + x_{cen}\\
        b_{y} \sin{\left( t \right)}
    \end{array} \right),
    \label{eq:InputSys}
\end{equation}
for $t=0, \tau, 2 \tau, \ldots,$ using the time-step $\tau = 0.01$. The resultant time-series of $\boldsymbol{u}(t)$ corresponds to a trajectory around a circle of radius $b=|b_{x}|=|b_{y}|$ and centered at $\left( x_{cen}, \, 0 \right)$. 

As in \cite{flynn2023seeingdouble}, for a given $x_{cen}$ we set $b_{x}=b_{y}=5$ in Eq.~\eqref{eq:InputSys} to generate a trajectory about the counter-clockwise circular orbit that we denote as $\mathcal{C}_{A}$ and points along this orbit are written as $\boldsymbol{u}_{\left(\mathcal{C}_{A}\right)}(t)$. For the corresponding $-x_{cen}$, we generate a trajectory about the clockwise circular orbit that we denote as $\mathcal{C}_{B}$, by setting $b_{x}=-5$ and $b_{y}=5$, points along this orbit are written as $\boldsymbol{u}_{\left(\mathcal{C}_{B}\right)}(t)$. 
By changing $x_{cen}$ the centres of these cycles are moved equidistantly along the line $y=0$. An overlapping region between $\mathcal{C}_{A}$ and $\mathcal{C}_{B}$ exists whenever $|x_{cen}| < b = 5$, i.e., $\mathcal{C}_{A} \cap \mathcal{C}_{B} \neq \emptyset$ $\forall \, |x_{cen}| < 5$. Furthermore, $\mathcal{C}_{A}$ and $\mathcal{C}_{B}$ are said to be `entirely/completely overlapping' when $x_{cen}=0$. In this extreme case, the only difference between $\mathcal{C}_{A}$ and $\mathcal{C}_{B}$ is the direction of rotation on both cycles.

The values of $\boldsymbol{u}_{\left(\mathcal{C}_{A}\right)}(t)$ and $\boldsymbol{u}_{\left(\mathcal{C}_{B}\right)}(t)$ are used as the input to the open-loop RC in Eq.\,\eqref{eq:ListenRes} for $0 \leq t \leq t_{\text{train}}$. The open-loop RC's response to these driving input signals are denoted as $\boldsymbol{r}_{\left(\mathcal{C}_{A}\right)}(t)$ and $\boldsymbol{r}_{\left(\mathcal{C}_{B}\right)}(t)$. Following the steps outlined in Sec.\,\ref{sssec:RCtraining}, the values of $\boldsymbol{r}_{\left(\mathcal{C}_{A}\right)}(t)$, $\boldsymbol{r}_{\left(\mathcal{C}_{B}\right)}(t)$, $\boldsymbol{u}_{\left(\mathcal{C}_{A}\right)}(t)$, and $\boldsymbol{u}_{\left(\mathcal{C}_{B}\right)}(t)$ for $t \in \left[ t_{\text{listen}}, t_{\text{train}} \right]$ are used to produce the corresponding training data matrices $\textbf{X}_{\left(\mathcal{C}_{A}\right)}, \, \textbf{X}_{\left(\mathcal{C}_{B}\right)}, \, \textbf{Y}_{\left(\mathcal{C}_{A}\right)},$ and $\textbf{Y}_{\left(\mathcal{C}_{B}\right)}$ as per Eqs.~\eqref{eq:Xmat} and \eqref{eq:Ymat} in order to compute $\textbf{W}_{out}$ in Eq.\,\eqref{eq:WoutRegression}. This $\textbf{W}_{out}$ is then used to create the closed-loop RC in Eq.\,\eqref{eq:PredRes}.

We say that this closed-loop RC achieves multifunctionality and solves the seeing double problem once it reconstructs a coexistence of $\mathcal{C}_{A}$ and $\mathcal{C}_{B}$. To do this the RC must construct a coexistence of two attractors, $\mathcal{S}_{A}$ and $\mathcal{S}_{B}$, that exist in $\mathbb{S}$ and resemble $\mathcal{C}_{A}$ and $\mathcal{C}_{B}$ when projected to $\mathbb{P}$ using $\hat{\boldsymbol{\psi}}\left(\cdot\right)$ in Eq.\,\eqref{eq:ReadoutFunction} with $\textbf{W}_{out}$ computed as mentioned above. As per the same convention used earlier, the projected dynamics of $\mathcal{S}_{A}$ and $\mathcal{S}_{B}$ are referred to as the reconstructed attractors and are denoted by $\hat{\mathcal{C}}_{A}$ and $\hat{\mathcal{C}}_{B}$. To reconstruct the dynamics of $\mathcal{C}_{A}$ or $\mathcal{C}_{B}$ using this multifunctional RC we initialise Eq.\,\eqref{eq:PredRes} with $\hat{\boldsymbol{r}}(0)=\boldsymbol{r}_{\left(\mathcal{C}_{A}\right)}(t_{\text{train}})$ or $\hat{\boldsymbol{r}}(0)=\boldsymbol{r}_{\left(\mathcal{C}_{B}\right)}(t_{\text{train}})$ or some known point in the basin of attraction of $\mathcal{S}_{A}$ or $\mathcal{S}_{B}$. The subsequent states of Eq.\,\eqref{eq:PredRes} when approaching $\mathcal{S}_{A}, \mathcal{S}_{B} \subset \mathbb{S}$ (i.e., $\hat{\mathcal{C}}_{A}, \hat{\mathcal{C}}_{B} \subset \mathbb{P}$) are written as $\hat{\boldsymbol{r}}_{\left(\mathcal{C}_{A}\right)}(t)$ and $\hat{\boldsymbol{r}}_{\left(\mathcal{C}_{B}\right)}(t)$.

\section{Results}

\subsection{Outline of experiments}\label{sec:Results_ExperOutline}

\begin{figure}
    \centering
    \includegraphics[width=.55\textwidth]{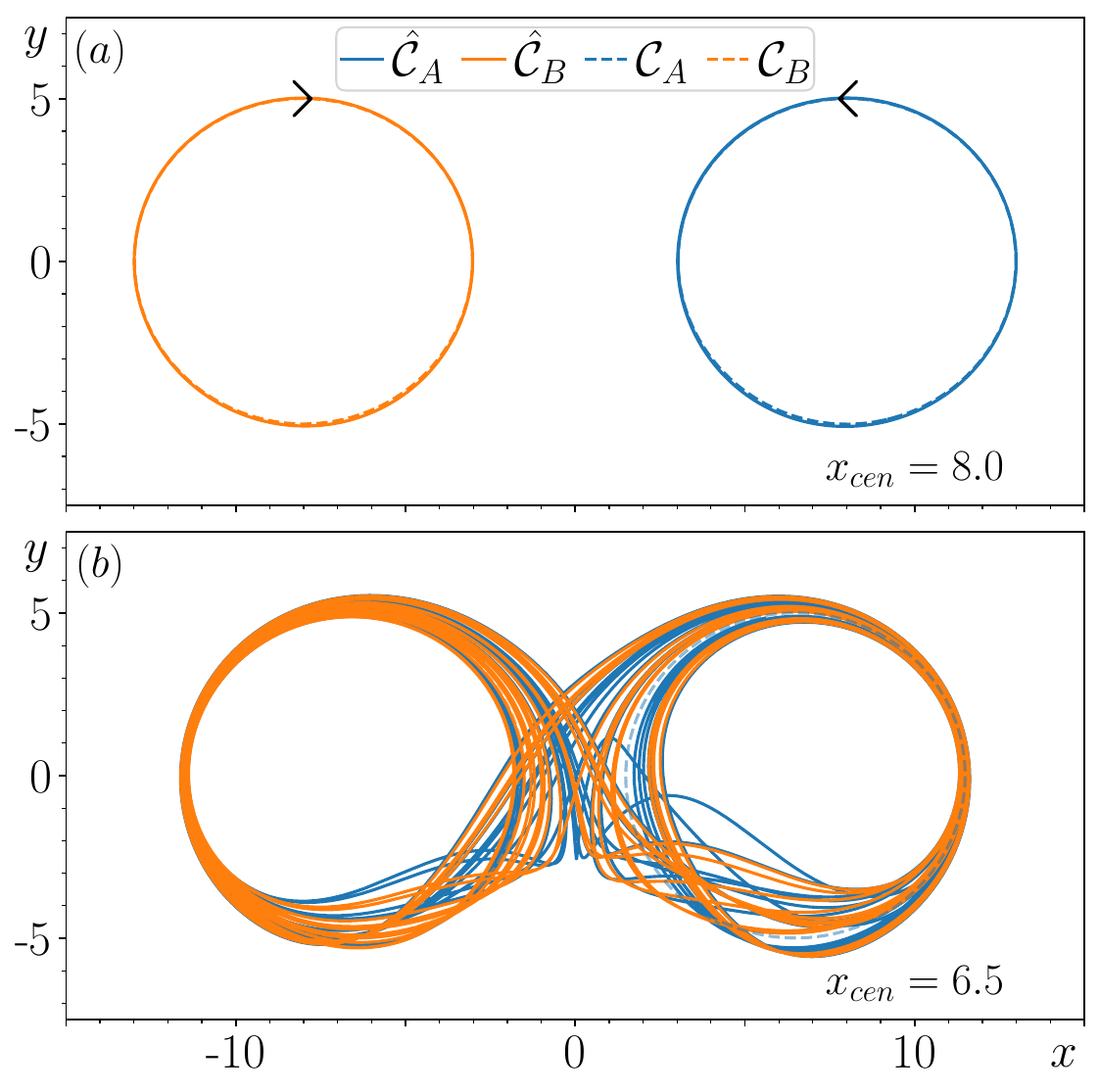}
    \caption{Illustrating the result of training the RC to reconstruct a coexistence of $\mathcal{C}_{A}$ and $\mathcal{C}_{B}$ when $\rho=0.2$ for $x_{cen}=8.0$ in panel (a) and for $x_{cen}=6.5$ in panel (b). Black arrows indicate direction of rotation on both orbits. Dynamics of the closed-loop RC are illustrated in solid curves, training data by dashed curves.}
    \label{fig:xytraj_examples_xcen_80_65_rho_02_}
\end{figure}

The main aim of this paper is to improve our current understanding of how metastable switching dynamics emerge in a RC that fails to achieve multifunctionality, Fig.\,\ref{fig:xytraj_examples_xcen_80_65_rho_02_} illustrates the particular phenomenon we are interested in studying. In panel (a) we show that when $\mathcal{C}_{A}$ and $\mathcal{C}_{B}$ are sufficiently far apart (when $x_{cen} = 8.0$) then for $\rho = 0.2$ the closed-loop RC achieves multifunctionality as $\hat{\mathcal{C}}_{A}$ and $\hat{\mathcal{C}}_{B}$ are more or less identical to $\mathcal{C}_{A}$ and $\mathcal{C}_{B}$. However panel (b) shows that when the same RC is trained with $x_{cen}=6.5$, when $\mathcal{C}_{A}$ and $\mathcal{C}_{B}$ are slightly closer together but do not overlap, then the closed-loop RC fails to achieve multifunctionality and instead its state switches between regions of $\mathbb{P}$ associated with $\mathcal{C}_{A}$ and $\mathcal{C}_{B}$. To investigate the origins of these switching dynamics we conduct the following experiments.

The results from the experiments reported on in this section consist of training the RC in Eq.~\eqref{eq:ListenRes} to solve the seeing double problem for $x_{cen} = 6.5,\, 5.0,\, 3.5,\,$ and $2.0$. To illustrate the differences in the closed-loop RC's (Eq.\eqref{eq:PredRes}) dynamics when trained at these values of $x_{cen}$, we chose a common $\rho$ value ($\rho = 0.7$) where multifunctionality is achieved for each $x_{cen}$. We then decrease $\rho$ in small steps of $0.001$ and track the changes in the dynamics of the reconstructed attractors, $\hat{\mathcal{C}}_{A}$ and $\hat{\mathcal{C}}_{B}$, by initialising the closed-loop RC with an initial condition corresponding to each attractor at the previous step and integrating the closed-loop RC forward in time up to $t=200$. If $\hat{\mathcal{C}}_{A}$ or $\hat{\mathcal{C}}_{B}$ can no longer be tracked, i.e. have become unstable or cease to exist, we continue decreasing $\rho$ and track the changes in the attractor that the state of the closed-loop RC subsequently approaches until $\rho = 0.1$. This method of attractor continuation enables us to investigate the origin of the switching dynamics we see in $\mathbb{P}$ at certain values of $\rho$ and $x_{cen}$. 

The result of this continuation procedure at each of the specified $x_{cen}$ values are shown in panel (e) of Figs.~\ref{fig:bifplot_xcen6.5_} - \ref{fig:bifplot_xcen2.0_} where we plot the local maxima of the reconstructed $x$ variable, denoted by $x_{m}$. In panels (a)-(d) of Figs.~\ref{fig:bifplot_xcen6.5_} - \ref{fig:bifplot_xcen2.0_} we illustrate some of the most significant changes in the closed-loop RC's dynamics at particular $\rho$ values, highlighting how the switching dynamics emerge. In Figs.~\ref{fig:bifplot_xcen6.5_} - \ref{fig:bifplot_xcen2.0_} the dashed blue and orange curves illustrate the location of $\mathcal{C}_{A}$ and $\mathcal{C}_{B}$, the corresponding solid curves describe the closed-loop RC's reconstruction of $\mathcal{C}_{A}$ and $\mathcal{C}_{B}$ (denoted by $\hat{\mathcal{C}}_{A}$ and $\hat{\mathcal{C}}_{B}$), and the blue and orange points are the corresponding $x_{m}$ values obtained from tracking the changes in the dynamics of $\hat{\mathcal{C}}_{A}$ and $\hat{\mathcal{C}}_{B}$ and the subsequent attractor that the closed-loop RC's state approaches when it fails to reconstruct $\hat{\mathcal{C}}_{A}$ or $\hat{\mathcal{C}}_{B}$. In circumstances where the closed-loop RC fails to reconstruct $\mathcal{C}_{A}$, $\mathcal{C}_{B}$, or produce switching dynamics between regions of $\mathbb{P}$ associated $\mathcal{C}_{A}$ and $\mathcal{C}_{B}$, the `untrained attractor' (an attractor that the closed-loop RC produces that was not present during the training) that the state of the closed-loop RC subsequently approaches is depicted using the colour specified in the associated plot legends. While there may be other untrained attractors present in $\mathbb{P}$ for $\rho < 0.1$, in order to maintain the focus of this paper (which is to explore the origins of the switching dynamics), we only track the changes in the attractors that the state of the closed-loop RC approaches when it fails to reconstruct $\mathcal{C}_{A}$, $\mathcal{C}_{B}$, or produce the switching dynamics. Saying that, we also initialise the state of the closed-loop RC from many random initial conditions at several different $\rho$ values when tracking the changes in $\hat{\mathcal{C}}_{A}$ and $\hat{\mathcal{C}}_{B}$ but we do not find any untrained attractors.

In panel (e) of Figs.~\ref{fig:bifplot_xcen6.5_} - \ref{fig:bifplot_xcen2.0_} the vertical dashed grey lines indicate the $\rho$ values that the corresponding dynamics in $\mathbb{P}$ are illustrated in the accompanying panels (a)-(d). The black arrows in panel (a) of Figs.~\ref{fig:bifplot_xcen6.5_} - \ref{fig:bifplot_xcen2.0_} indicate $\mathcal{C}_{A}$ and $\mathcal{C}_{B}$'s direction of rotation. These illustrations are generated by training the RC in Eq.\eqref{eq:ListenRes} at the specified $\rho$ values and initialising the closed-loop RC with $\boldsymbol{r}_{\left(\mathcal{C}_{A}\right)}(t_{\text{train}})$ and $\boldsymbol{r}_{\left(\mathcal{C}_{B}\right)}(t_{\text{train}})$ (the last point in the training data corresponding to $\mathcal{C}_{A}$ and $\mathcal{C}_{B}$), i.e., following the description in Sec.~\ref{sssec:RCtraining}. We do this in order to observe whether there exists any transient behaviour associated with $\mathcal{C}_{A}$ or $\mathcal{C}_{B}$ when the RC fails to reconstruct these attractors.

Furthermore, we also conducted this analysis for the case of $x_{cen} = 8.0$ but we choose not show these results as we found no switching dynamics nor significant changes in the closed-loop RC's dynamics for changes in $\rho$. For $x_{cen} = 8.0$, the closed-loop RC achieves multifunctionality with nearly perfect reconstruction of $\mathcal{C}_{A}$ and $\mathcal{C}_{B}$ for the range of $\rho$ values that were investigated. 

The results of the continuation analysis for each of the selected $x_{cen}$ values are outlined in Secs.~\ref{sec:res_xcen65}-\ref{sec:res_xcen20}. In Sec.~\ref{ssec:closer_inspections} we examine the residence and escape times associated with the switching and transient dynamics observed in Secs.~\ref{sec:res_xcen65}-\ref{sec:res_xcen20}.

\subsection{Continuation analysis for $x_{cen}=6.5$}\label{sec:res_xcen65}

\begin{figure}
    \centering
    \includegraphics[width=.85\textwidth]{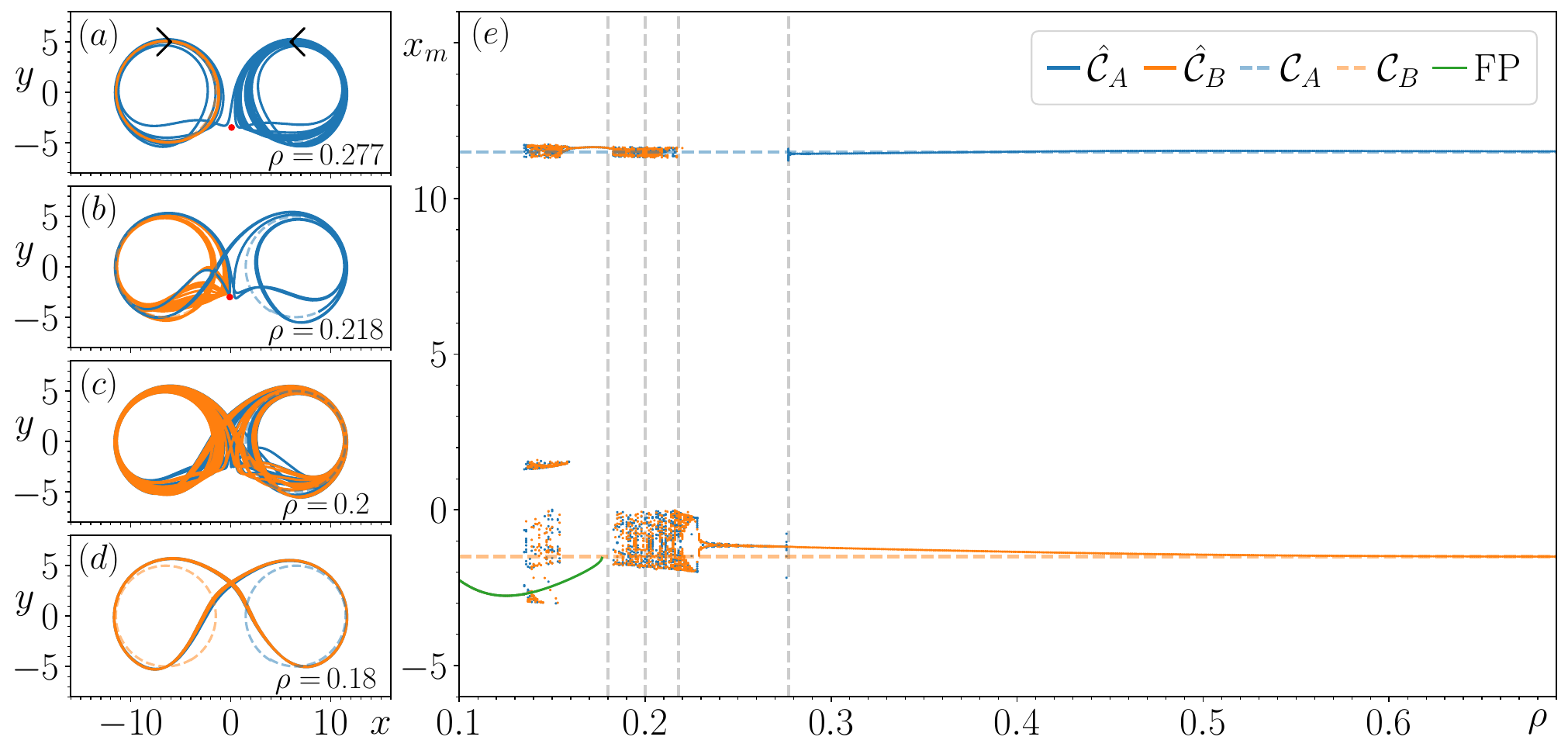}
    \caption{Illustrating the result of tracking the changes in $\hat{\mathcal{C}}_{A}$ and $\hat{\mathcal{C}}_{B}$ with respect to changes in $\rho$ for $x_{cen}=6.5$. Panel (e) describes how the local maxima of the corresponding attractors that are tracked, $x_{m}$, changes with respect to $\rho$. Panels (a)-(d) highlight some of the most significant changes in the dynamics of $\hat{\mathcal{C}}_{A}$ and $\hat{\mathcal{C}}_{B}$ at certain values of $\rho$ from the perspective of $\mathbb{P}$, the prediction state space.}
    \label{fig:bifplot_xcen6.5_}
\end{figure}

Fig.~\ref{fig:bifplot_xcen6.5_}~(e) shows that for $x_{cen}=6.5$ there are no significant changes in the closed-loop RC's ability to reconstruct a coexistence of $\mathcal{C}_{A}$ and $\mathcal{C}_{B}$ until $\rho = 0.277$ where $\hat{\mathcal{C}}_{A}$ becomes unstable. In Fig.~\ref{fig:bifplot_xcen6.5_}~(a) we take a closer look at this behaviour, we see that the state of the closed-loop RC follows a significantly long chaotic transient when initialised with the associated $\boldsymbol{r}_{\left(\mathcal{C}_{A}\right)}(t_{\text{train}})$ before eventually approaching $\hat{\mathcal{C}}_{B}$, which is the only stable attractor present in $\mathbb{P}$ (confirmed by initialising the closed-loop RC from many random initial conditions). Based on the structure of this transient and the evidence of a saddle present at $(x,y) \approx (0.1, -3.5)$ (indicated by the red point in Fig.~\ref{fig:bifplot_xcen6.5_}~(a)) there is evidence that $\hat{\mathcal{C}}_{A}$ becomes unstable by colliding with this saddle. 

There is stronger evidence to support the above claim in Fig.~\ref{fig:bifplot_xcen6.5_}~(b) as we see for $\rho=0.218$ the structure of the transient is highly influenced by this saddle (now located at $(x,y) \approx (0, -3)$ in $\mathbb{P}$ as indicated by the red point in Fig.~\ref{fig:bifplot_xcen6.5_}~(b)) whose unstable directions appear to point predominately along the x-axis and stable directions point along the y-axis in $\mathbb{P}$. While we see here that the state of the closed-loop RC takes a significantly shorter amount of time to escape this transient, what is particularly interesting about this transient is that the state of the closed-loop RC follows a path that encircles $\hat{\mathcal{C}}_{B}$ and switches back to the portion of $\mathbb{P}$ associated with $\mathcal{C}_{A}$ before switching back to and remaining on $\hat{\mathcal{C}}_{B}$ for all future time. In fact, $\hat{\mathcal{C}}_{B}$ is the only stable attractor present in $\mathbb{P}$ for $\rho \in [0.218, 0.277]$ as further indicated by the dashed blue horizontal line associated with the reconstruction of $\mathcal{C}_{A}$ that is visible within this range of $\rho$ values. It is also during this range of $\rho$ values that $\hat{\mathcal{C}}_{B}$ becomes chaotic through what appears to be a period-doubling bifurcation that is quickly followed by a torus bifurcation. Furthermore, it appears that the closed-loop RC's trajectory on $\hat{\mathcal{C}}_{B}$ for $\rho = 0.218$ comes arbitrarily close to the saddle, indicating that a similar fate to $\hat{\mathcal{C}}_{A}$ awaits $\hat{\mathcal{C}}_{B}$ at smaller $\rho$ values.

Fig.~\ref{fig:bifplot_xcen6.5_}~(e) illustrates that when tracking the changes in $\hat{\mathcal{C}}_{B}$ for decreasing $\rho$ further there are multiple values of $x_{m}$ obtained that surround both blue and orange dashed horizontal lines. This indicates the emergence of the switching dynamics between regions of $\mathbb{P}$ where the previously stable $\hat{\mathcal{C}}_{A}$ and $\hat{\mathcal{C}}_{B}$ existed. These switching dynamics are found to occur for $\rho \in [0.135, 0.217]$ and there are a variety of different switching patterns exhibited. For instance, in Fig.~\ref{fig:bifplot_xcen6.5_}~(c) we show that these switching dynamics resemble a Lorenz-like chaotic attractor for $\rho = 0.2$, whereas in Fig.~\ref{fig:bifplot_xcen6.5_}~(d) a periodic switching pattern appears in the form of a limit cycle, which emerges from the chaotic attractor. In Sec.~\ref{ssec:closer_inspections} we take a closer looks at the long-term dynamics of the chaotic attractor shown in Fig.~\ref{fig:bifplot_xcen6.5_}~(c). Fig.~\ref{fig:bifplot_xcen6.5_}~(e) shows that this periodic switching pattern returns to an aperiodic switching pattern at $\rho = 0.155$ indicated by the three clusters of $x_{m}$ values.

The switching dynamics come to an end at $\rho=0.135$ and the state of the closed-loop RC subsequently approaches a stable fixed point (FP), indicated by the sequence of green points, and we continue to track the changes in this FP until $\rho=0.1$. We find a small range of $\rho$ values where the FP coexists with the switching patterns by tracking the changes in this FP for increasing $\rho$ values until it becomes unstable at $\rho = 0.178$. At this point the state of the closed-loop RC returns to the limit cycle associated with the periodic switching pattern mentioned in the paragraph above.

\subsection{Continuation analysis for $x_{cen}=5.0$}\label{sec:res_xcen50}

\begin{figure}
    \centering
    \includegraphics[width=.85\textwidth]{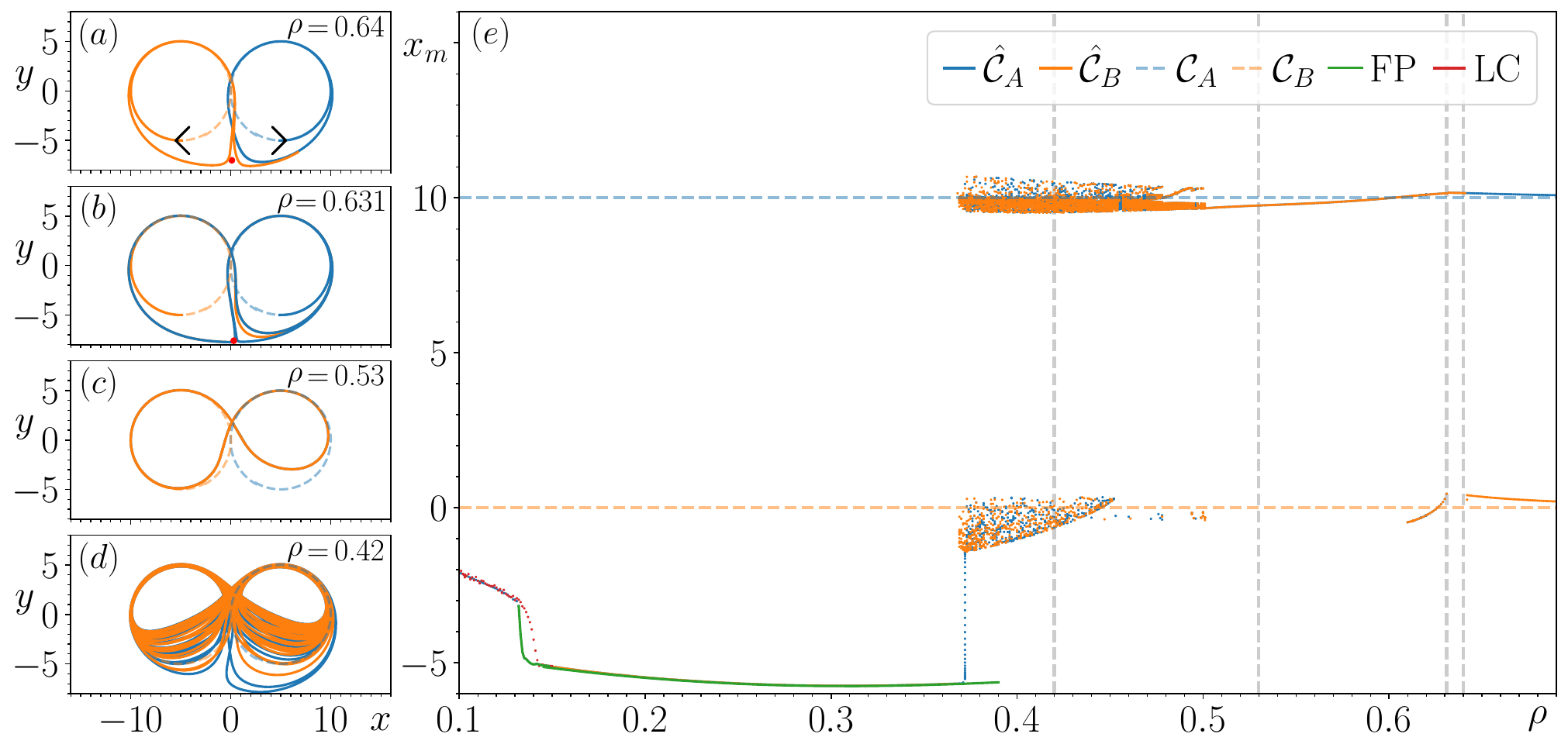}
    \caption{Illustrating the result of tracking the changes in $\hat{\mathcal{C}}_{A}$ and $\hat{\mathcal{C}}_{B}$ with respect to changes in $\rho$ for $x_{cen}=5.0$. Panel (e) describes how the local maxima of the corresponding attractors that are tracked, $x_{m}$, changes with respect to $\rho$. Panels (a)-(d) highlight some of the most significant changes in the dynamics of $\hat{\mathcal{C}}_{A}$ and $\hat{\mathcal{C}}_{B}$ at certain values of $\rho$ from the perspective of $\mathbb{P}$, the prediction state space.}
    \label{fig:bifplot_xcen5.0_}
\end{figure}

For $x_{cen} = 5.0$ we find that by moving $\mathcal{C}_{A}$ and $\mathcal{C}_{B}$ closer together so that they touch at $(x,y)=(0,0)$, the switching dynamics begin at much larger $\rho$ values, persist for a greater range of $\rho$ values, and the switching patterns are first found to occur periodically before becoming chaotic. Fig.~\ref{fig:bifplot_xcen5.0_}~(e) shows that as $\rho$ is decreased, $\hat{\mathcal{C}}_{B}$ becomes unstable at $\rho = 0.64$. Fig.~\ref{fig:bifplot_xcen5.0_}~(a) illustrates the transient dynamics of the closed-loop RC when initialised with the associated $\boldsymbol{r}_{\left(\mathcal{C}_{B}\right)}(t_{\text{train}})$. Here we see that the state of the closed-loop RC follows one loop around the dashed orange circle ($\mathcal{C}_{B}$) before diverging away to approach and then remain on the slightly oval-shaped $\hat{\mathcal{C}}_{A}$. It is from these unstable directions of flow along the x-axis and stable flow along the y-axis that the nature of this transient also provides us with some evidence that there is a saddle located at $\mathcal{C}_{B}$ ($(x,y) \approx (0.1, -7)$) as indicated by the red point in Fig.~\ref{fig:bifplot_xcen5.0_}~(a).

Fig.~\ref{fig:bifplot_xcen5.0_}~(b) provides us with further information about this saddle (now located at $(x,y) \approx (0.3, -7.6)$ as indicated by the red point in Fig.~\ref{fig:bifplot_xcen5.0_}~(b)) as it appears that $\hat{\mathcal{C}}_{A}$ has become unstable at $\rho=0.631$ by colliding with the saddle. Moreover, it is through this second collision that a new limit cycle is created that produces the switching dynamics nearby the point at which $\mathcal{C}_{A}$ and $\mathcal{C}_{B}$ touch in $\mathbb{P}$. This limit cycle consists of two weakly attracting connected regions of flow around $\mathcal{C}_{A}$ and $\mathcal{C}_{B}$. Taking a closer look at the transient dynamics exhibited by the closed-loop RC when initialised with $\boldsymbol{r}_{\left(\mathcal{C}_{A}\right)}(t_{\text{train}})$, we see that its state comes arbitrarily close to the saddle before completing one loop around the dashed blue circle associated with $\mathcal{C}_{A}$, however on the second loop the trajectory diverges away from $\mathcal{C}_{A}$ nearby the saddle and then approaches and subsequently remains on the new larger limit cycle that consists of loops around regions of $\mathbb{P}$ associated with $\mathcal{C}_{A}$ and $\mathcal{C}_{B}$. Initially there are two values obtained for $x_{m}$, the local maxima associated with $\mathcal{C}_{A}$ and a point nearby the saddle, the small branch of points nearby the dashed orange horizontal line seen in Fig.~\ref{fig:bifplot_xcen5.0_}~(e). Correspondingly, the sharp turning point on the new limit cycle nearby the saddle point shown in Fig.~\ref{fig:bifplot_xcen5.0_}~(b) does not persist for many subsequent $\rho$ values as the limit cycle starts to resemble a figure of 8 shape in $\mathbb{P}$ like the example shown in Fig.~\ref{fig:bifplot_xcen5.0_}~(c) for $\rho = 0.53$. 

As shown in Fig.~\ref{fig:bifplot_xcen5.0_}~(e), additional values of $x_{m}$ are found for $\rho = 0.5$ as the limit cycle transitions to a chaotic attractor. In Fig.~\ref{fig:bifplot_xcen5.0_}~(d) we provide an example of the aperiodic switching patterns exhibited by this chaotic attractor for $\rho = 0.42$. In Sec.~\ref{ssec:closer_inspections} we take a closer look at the long-term dynamics of this chaotic attractor. We are unable to track the changes in this chaotic attractor for $\rho < 0.37$ and the state of the closed-loop RC subsequently approaches a FP, whose behaviour with respect to changes in $\rho$ is described the green branch of points in Fig.~\ref{fig:bifplot_xcen5.0_}~(e). There is a relatively small range of $\rho$ values where this FP coexists with the chaotic attractor associated with the aperiodic switching patterns for $\rho \in [0.37, 0.39]$. When tracking the changes in this FP for decreasing $\rho$ we find a smaller range of $\rho$ values where this FP coexists with a different period-1 limit cycle, whose corresponding $x_{m}$ is described by the branch of red points in Fig.~\ref{fig:bifplot_xcen5.0_}~(e).

\subsection{Continuation analysis for $x_{cen}=3.5$}\label{sec:res_xcen35}

\begin{figure}
    \centering
    \includegraphics[width=.85\textwidth]{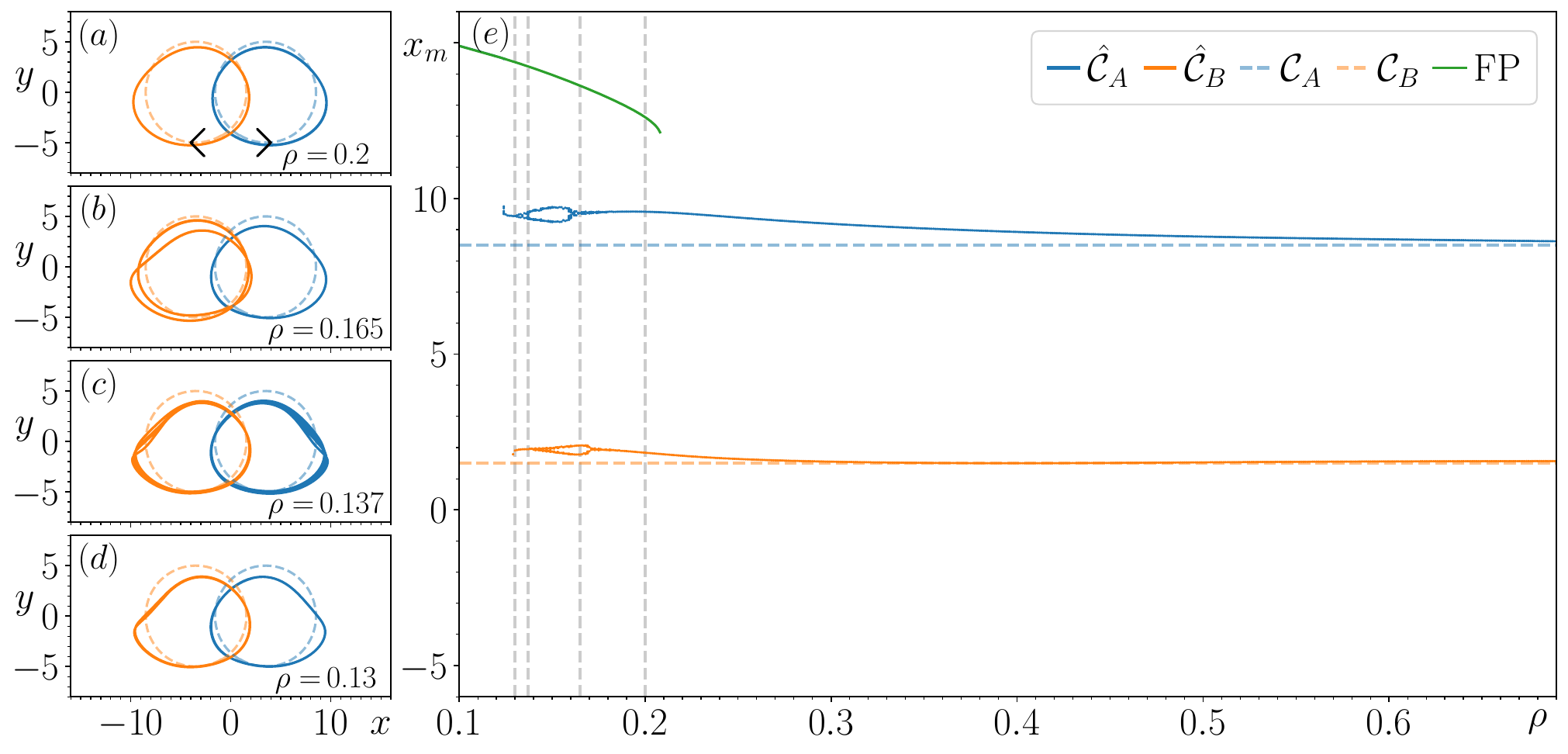}
    \caption{Illustrating the result of tracking the changes in $\hat{\mathcal{C}}_{A}$ and $\hat{\mathcal{C}}_{B}$ with respect to changes in $\rho$ for $x_{cen}=3.5$. Panel (e) describes how the local maxima of the corresponding attractors that are tracked, $x_{m}$, changes with respect to $\rho$. Panels (a)-(d) highlight some of the most significant changes in the dynamics of $\hat{\mathcal{C}}_{A}$ and $\hat{\mathcal{C}}_{B}$ at certain values of $\rho$ from the perspective of $\mathbb{P}$, the prediction state space.}
    \label{fig:bifplot_xcen3.5_}
\end{figure}

When decreasing $x_{cen}$ to $3.5$ we find that by moving $\mathcal{C}_{A}$ and $\mathcal{C}_{B}$ closer together so that they overlap and share two common points, this improves the performance of the RC as it is achieves multifunctionality for a much larger range of $\rho$ values and does not produce any switching dynamics. Fig.~\ref{fig:bifplot_xcen3.5_}~(e) shows that by tracking the changes in $\hat{\mathcal{C}}_{A}$ and $\hat{\mathcal{C}}_{B}$ for decreasing $\rho$ there is a growing difference between the obtained values for $x_{m}$ and the corresponding true values with respect to $\mathcal{C}_{A}$ and $\mathcal{C}_{B}$. Fig.~\ref{fig:bifplot_xcen3.5_}~(a) provides further insight into how $\hat{\mathcal{C}}_{A}$ and $\hat{\mathcal{C}}_{B}$ deform as $\rho$ decreases. Fig.~\ref{fig:bifplot_xcen3.5_}~(b) shows how $\hat{\mathcal{C}}_{A}$ and $\hat{\mathcal{C}}_{B}$ increasingly lose their resemblance to $\mathcal{C}_{A}$ and $\mathcal{C}_{B}$ with $\hat{\mathcal{C}}_{B}$ having undergone a period-doubling bifurcation as $\rho$ is decreased to $\rho = 0.165$. Fig.~\ref{fig:bifplot_xcen3.5_}~(c) illustrates that for $\rho = 0.137$ both $\hat{\mathcal{C}}_{A}$ and $\hat{\mathcal{C}}_{B}$ display aperiodic dynamics (indicated by the increased thickness of the corresponding blue and orange curves). Fig.~\ref{fig:bifplot_xcen3.5_}~(e) shows that as $\rho$ is decreased further, $\hat{\mathcal{C}}_{B}$ becomes unstable at $\rho \approx 0.129$ and the state of the closed-loop RC subsequently approaches $\hat{\mathcal{C}}_{A}$, at $\rho \approx 0.124$ we find that $\hat{\mathcal{C}}_{A}$ becomes unstable and the state of the closed-loop RC subsequently approaches the FP described by the branch of green points. By tracking the changes in this FP for increasing $\rho$ we find that this FP coexists with $\hat{\mathcal{C}}_{A}$ and $\hat{\mathcal{C}}_{B}$ until it becomes unstable at $\rho \approx 0.21$ and the state of the closed-loop RC returns to $\hat{\mathcal{C}}_{A}$. Fig.~\ref{fig:bifplot_xcen3.5_}~(d) illustrates that prior to $\hat{\mathcal{C}}_{A}$ and $\hat{\mathcal{C}}_{B}$ becoming unstable these attractors are no longer chaotic and have returned to period-1 limit cycles.

\subsection{Continuation analysis for $x_{cen}=2.0$}\label{sec:res_xcen20}

\begin{figure}
    \centering
    \includegraphics[width=.85\textwidth]{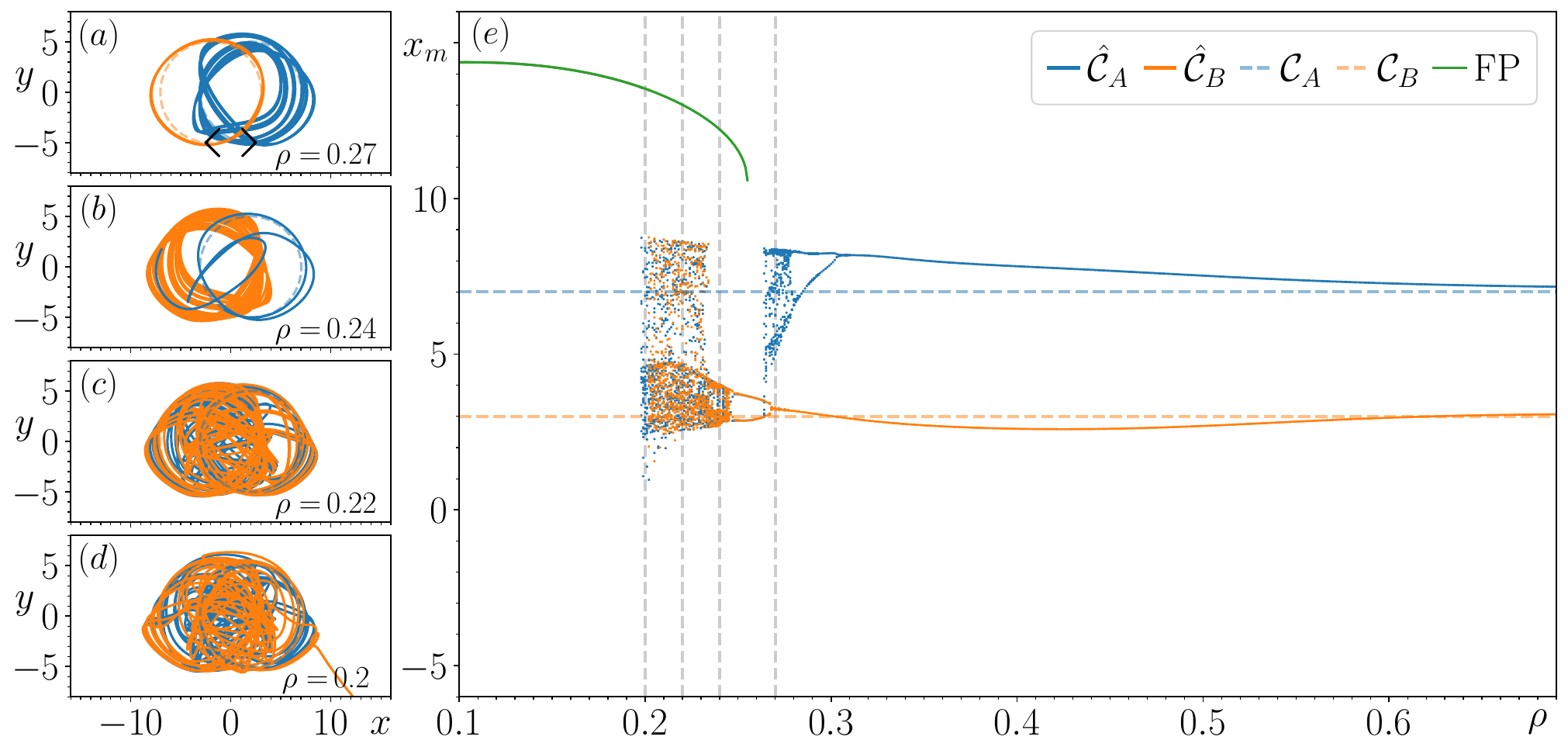}
    \caption{Illustrating the result of tracking the changes in $\hat{\mathcal{C}}_{A}$ and $\hat{\mathcal{C}}_{B}$ with respect to changes in $\rho$ for $x_{cen}=2.0$. Panel (e) describes how the local maxima of the corresponding attractors that are tracked, $x_{m}$, changes with respect to $\rho$. Panels (a)-(d) highlight some of the most significant changes in the dynamics of $\hat{\mathcal{C}}_{A}$ and $\hat{\mathcal{C}}_{B}$ at certain values of $\rho$ from the perspective of $\mathbb{P}$, the prediction state space.}
    \label{fig:bifplot_xcen2.0_}
\end{figure}

When decreasing $x_{cen}$ to $2.0$, we find that by increasing the amount of overlap between $\mathcal{C}_{A}$ and $\mathcal{C}_{B}$, the closed-loop RC produces switching dynamics within a relatively small range of $\rho$ values in a similar fashion to those found for $x_{cen} = 6.5$ and $5.0$. Fig.~\ref{fig:bifplot_xcen2.0_}~(e) shows that as $\rho$ decreases there is a increasingly large offset between the values of $x_{m}$ for $\hat{\mathcal{C}}_{A}$ and $\mathcal{C}_{A}$ and a small but noticeable difference between the values of $x_{m}$ for $\hat{\mathcal{C}}_{B}$ and $\mathcal{C}_{B}$. 

Fig.~\ref{fig:bifplot_xcen2.0_}~(e) shows that as $\rho$ is decreased from $\rho \approx 0.3$, $\hat{\mathcal{C}}_{A}$ undergoes a sequence of period-doubling bifurcations which results in $\hat{\mathcal{C}}_{A}$ transitions to a chaotic attractor. In Fig.~\ref{fig:bifplot_xcen2.0_}~(a) we illustrate the dynamics of the closed-loop RC for $\rho = 0.27$ which shows the coexistence of the chaotic $\hat{\mathcal{C}}_{A}$ and periodic $\hat{\mathcal{C}}_{B}$ which closely resembles the dynamics of $\mathcal{C}_{B}$. Fig.~\ref{fig:bifplot_xcen2.0_}~(e) shows that by decreasing $\rho$ further, $\hat{\mathcal{C}}_{B}$ undergoes a period-doubling bifurcation starting from $\rho \approx 0.269$. $\hat{\mathcal{C}}_{A}$ becomes unstable at $\rho \approx 0.261$ and after a bout of transient dynamics the state of the closed-loop RC subsequently approaches the period-2 $\hat{\mathcal{C}}_{B}$. We then continue to track the changes in $\hat{\mathcal{C}}_{B}$ which also becomes chaotic at $\rho \approx 0.245$. In Fig.~\ref{fig:bifplot_xcen2.0_}~(b) we illustrate the chaotic dynamics of $\hat{\mathcal{C}}_{B}$ and the relatively short duration of transient dynamics exhibited by the closed-loop RC when initialised from $\boldsymbol{r}_{\left(\mathcal{C}_{A}\right)}(t_{\text{train}})$. This transient completes one loop around the region of $\mathbb{P}$ associated with $\mathcal{C}_{A}$ however on its second loop the state of the closed-loop RC approaches the point $(x,y) \approx (-4.5,-4.5)$ where it subsequently reverses along its trajectory to this point and then approaches the chaotic $\hat{\mathcal{C}}_{B}$, remaining on $\hat{\mathcal{C}}_{B}$ thereafter.

The densely populated range of $x_{m}$ values which spans across both dashed horizontal lines associated with $\mathcal{C}_{A}$ and $\mathcal{C}_{B}$ in Fig.~\ref{fig:bifplot_xcen2.0_}~(e) shows that the switching dynamics emerge at $\rho \approx 0.23$. For $\rho = 0.22$, in Fig.~\ref{fig:bifplot_xcen2.0_}~(c) we illustrate the dynamics of the large chaotic attractor that is born at $\rho \approx 0.23$ and trajectories on this attractor resemble aperiodic switching dynamics between regions of $\mathbb{P}$ associated with $\mathcal{C}_{A}$ and $\mathcal{C}_{B}$. For $\rho=0.2$, Fig.~\ref{fig:bifplot_xcen2.0_}~(d) illustrates that when the closed-loop RC is initialised with $\boldsymbol{r}_{\left(\mathcal{C}_{B}\right)}(t_{\text{train}})$ its state follows a chaotic transient before approaching a FP located at $(x,y) \approx (13.5,-9.5)$ which is just outside the portion of $\mathbb{P}$ shown here. 

\subsection{Closer inspection of switching dynamics at $x_{cen}=6.5$ and $5.0$}\label{ssec:closer_inspections}

In this section we aim to shed further light on the nature of the switching dynamics discussed so far. We consider the two examples of $x_{cen}=6.5$ and $\rho = 0.2$ which we refer to as Case 1, and $x_{cen}=5.0$ and $\rho = 0.42$ which we refer to as Case 2. We generate a much longer trajectory on these chaotic attractors in order to determine the distribution of residence times that the state of the closed-loop RC spends in the respective $\mathcal{C}_{A}$ and $\mathcal{C}_{B}$ regions of $\mathbb{P}$. When the state of the closed-loop RC is in the region of $\mathbb{P}$ associated with $\mathcal{C}_{A}$ we consider the system to be in a metastable state denoted as $\tilde{\mathcal{C}}_{A}$, and similarly for $\mathcal{C}_{B}$ we consider to system to be in a different metastable state denoted as $\tilde{\mathcal{C}}_{B}$. 

\subsubsection{Algorithm to detect transitions}

In order to identify when the state of the closed-loop RC is in $\tilde{\mathcal{C}}_{A}$ or $\tilde{\mathcal{C}}_{B}$ we construct a relatively simple algorithm based on the concept of a `non-ideal relay' \cite{Pokrovskii12systemswHysteresis}. We use this algorithm to detect transition times from $\tilde{\mathcal{C}}_{A}$ to $\tilde{\mathcal{C}}_{B}$ and vice-versa. The non-ideal relay aspect of the algorithm involves choosing two threshold values, $\alpha$ and $\beta$, where we say that the closed-loop RC is in $\tilde{\mathcal{C}}_{A}$ once its state crosses $\beta$ and remains below $\alpha$, and is in $\tilde{\mathcal{C}}_{B}$ once its state crosses $\alpha$ and remains above $\beta$. The benefit of using these two thresholds as opposed to one threshold is that it allows us to improve our estimate of when the system is in a particular state by reducing the effect of false alarm scenarios where, for instance, the state of the closed-loop RC is in $\tilde{\mathcal{C}}_{A}$ but briefly dips below the single threshold and does not spend any significant amount of time in the portion of $\mathbb{P}$ associated with $\mathcal{C}_{B}$.

\begin{figure}
    \centering
    \includegraphics[width=.75\textwidth]{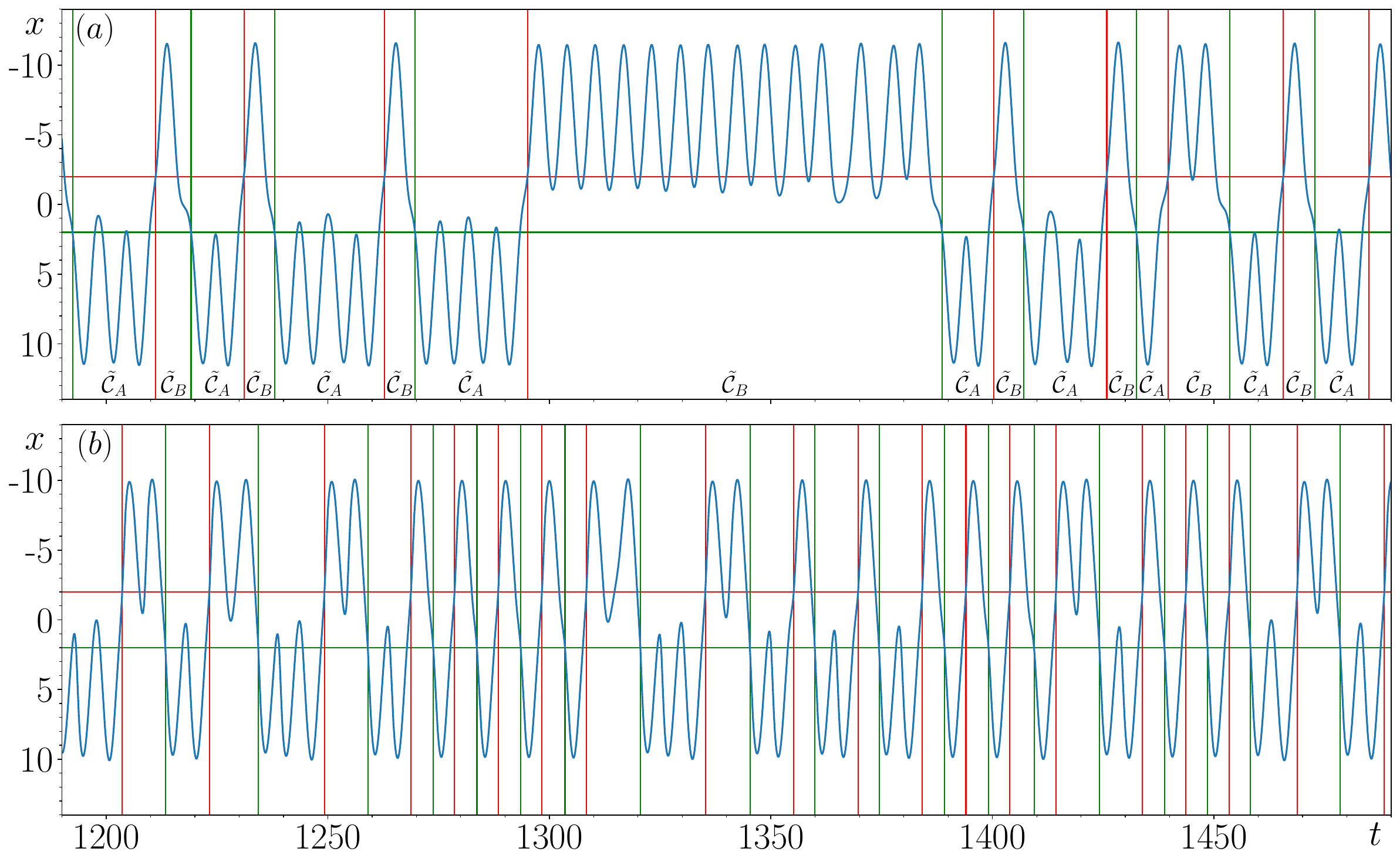}
    \caption{Illustrating how the transition times between $\tilde{\mathcal{C}}_{A}$ and $\tilde{\mathcal{C}}_{B}$ are obtained for Case 1 in (a) and Case 2 in (b).}
    \label{fig:xtraj_examples_CI_xcen_65_50_rho_020_042_}
\end{figure}

The result of using this algorithm to detect transitions from $\tilde{\mathcal{C}}_{A}$ to $\tilde{\mathcal{C}}_{B}$ and vice-versa in Case 1 is illustrated in Fig.~\ref{fig:xtraj_examples_CI_xcen_65_50_rho_020_042_}~(a) and for Case 2 is illustrated in Fig.~\ref{fig:xtraj_examples_CI_xcen_65_50_rho_020_042_}~(b) where we set $\alpha = - 2$ and $\beta = 2$. The green and red horizontal lines in Fig.~\ref{fig:xtraj_examples_CI_xcen_65_50_rho_020_042_} are used to illustrate these threshold values. The green and red vertical lines shown here correspond to the detected transitions times where the state of the closed-loop RC first enters $\tilde{\mathcal{C}}_{A}$ and $\tilde{\mathcal{C}}_{B}$ respectively. The residence times in $\tilde{\mathcal{C}}_{A}$ and $\tilde{\mathcal{C}}_{B}$ are then calculated based on these transition times.

The benefit of using this double threshold algorithm is made clear in Fig.~\ref{fig:xtraj_examples_CI_xcen_65_50_rho_020_042_}~(a), if instead a single threshold of $0$ were used then when the state of the RC crosses zero without switching from one metastable state to the other, like at $t \approx 1200, 1230, 1245, 1315, 1415, 1465, 1485$, then all of these crossings would be considered as transitions which is evidently false.

Furthermore, what is also evident from Fig.~\ref{fig:xtraj_examples_CI_xcen_65_50_rho_020_042_} is that there are at least three distinct types of switching pattern present where the state of the closed-loop RC can rapidly switch between $\tilde{\mathcal{C}}_{A}$ and $\tilde{\mathcal{C}}_{B}$ or spend a particular amount of time in $\tilde{\mathcal{C}}_{A}$ and $\tilde{\mathcal{C}}_{B}$ before switching. 

\subsubsection{Residence times}

\begin{figure}
    \centering
    \includegraphics[width=\textwidth]{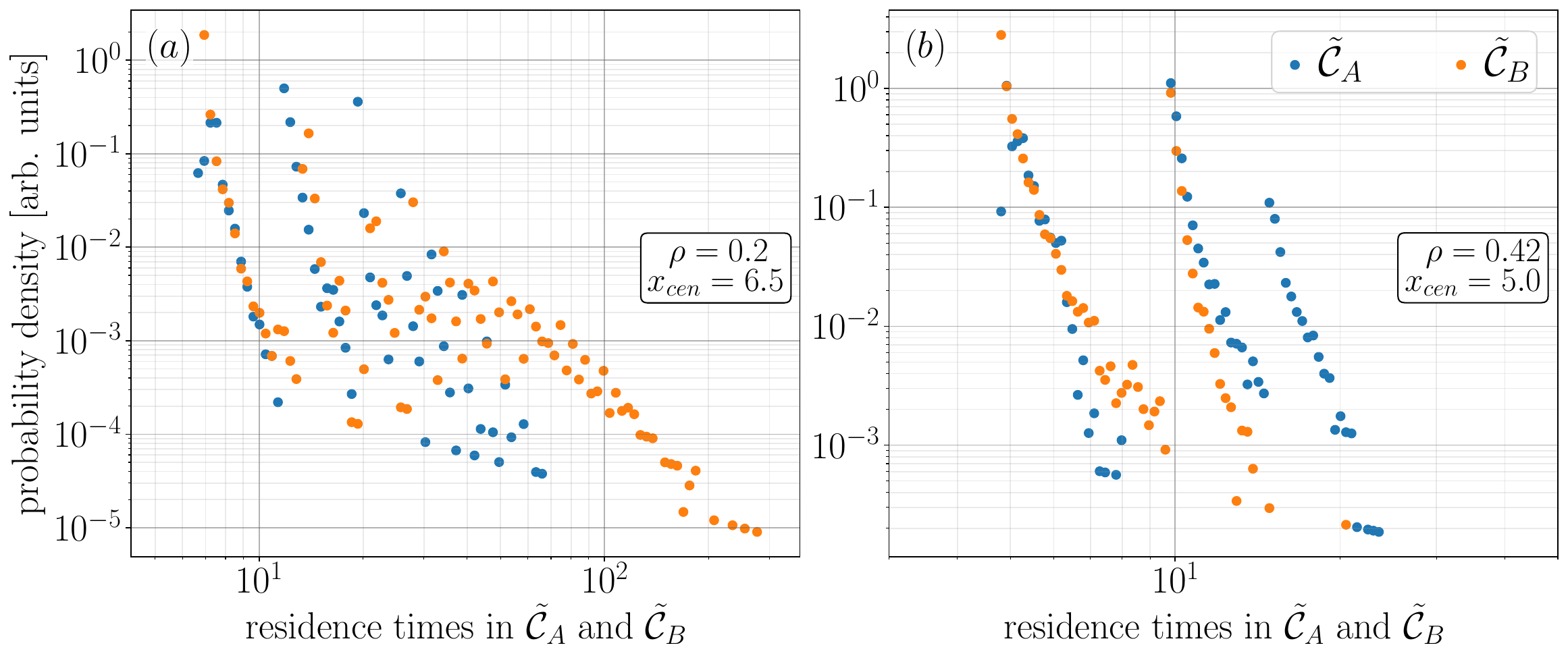}
    \caption{Probability density of residence times in $\tilde{\mathcal{C}}_{A}$ and $\tilde{\mathcal{C}}_{B}$ for Case 1 in (a) and Case 2 in (b).}
    \label{fig:ResidenceTimes_xcen_50_65_rho_042_02}
\end{figure}

In order to construct a reasonably well distributed sample of residence times in $\tilde{\mathcal{C}}_{A}$ and $\tilde{\mathcal{C}}_{B}$, we generate 10,000 examples of switchings between $\tilde{\mathcal{C}}_{A}$ and $\tilde{\mathcal{C}}_{B}$. To do this we needed to integrate the closed-loop RC forward in time up to $t \approx 300,000$ for Case 1 and up to $t \approx 157,000$ for Case 2. This tells us there are nearly twice as many switchings in a given duration of time for Case 2 in comparison to Case 1. From this sample of 10,000 switchings we found that for Case 1, the maximum and minimum residence time (in the arbitrary units of $t$) in $\tilde{\mathcal{C}}_{A}$ were $\approx 65$ and $6.5$ and for $\tilde{\mathcal{C}}_{B}$ were $\approx 270$ and $6.7$. For Case 2, the maximum and minimum residence time in $\tilde{\mathcal{C}}_{A}$ were $\approx 23.5$ and $4.8$ and for $\tilde{\mathcal{C}}_{B}$ were $\approx 20.3$ and $4.5$. We then compute the probability density of these residence times by generating a histogram of residence times with 100 bins chosen from numbers spaced evenly on a log scale with limits set to the max and min values specified above. The resulting probability density of these residence times for Case 1 is shown in Fig.~\ref{fig:ResidenceTimes_xcen_50_65_rho_042_02}~(a) and for Case 2 is shown in Fig.~\ref{fig:ResidenceTimes_xcen_50_65_rho_042_02}~(b).

What is most striking about the results shown in Fig.~\ref{fig:ResidenceTimes_xcen_50_65_rho_042_02} is that there is no single branch of exponentially distributed points, instead for both Case 1 and 2 the probability density of the residence times in $\tilde{\mathcal{C}}_{A}$ and $\tilde{\mathcal{C}}_{B}$ are organised into a number of branches of exponentially distributed points. 

We first outline the results shown in Fig.~\ref{fig:ResidenceTimes_xcen_50_65_rho_042_02}~(b) for Case 2 as it is more straightforward to discuss. The probability density of residence times in $\tilde{\mathcal{C}}_{A}$ are organised into three branches of exponentially distributed points and two branches of exponentially distributed points for $\tilde{\mathcal{C}}_{B}$. From further investigation we find that the points on these different branches correspond to scenarios where the state of the closed-loop RC follows either one or two loops (or partial loops) around $\mathcal{C}_{B}$ before switching to $\mathcal{C}_{A}$ and can follow up to three loops (or partial loops) about $\mathcal{C}_{A}$ before switching to $\mathcal{C}_{B}$. By partial loops we mean that the state of the closed-loop RC may switch from $\tilde{\mathcal{C}}_{A}$ to $\tilde{\mathcal{C}}_{B}$ without completing a full loop around $\mathcal{C}_{A}$. Furthermore, from the dynamics of the chaotic attractor which produces these switching dynamics illustrated in Fig.~\ref{fig:bifplot_xcen5.0_}~(d), it is reasonable to have anticipated the exponential distribution of point on these branches shown in Fig.~\ref{fig:ResidenceTimes_xcen_50_65_rho_042_02}~(b). It is also reasonable to have anticipated that the state of the closed-loop RC spends slightly longer amounts of time in $\tilde{\mathcal{C}}_{A}$ than $\tilde{\mathcal{C}}_{B}$ since $\hat{\mathcal{C}}_{B}$ becomes unstable before $\hat{\mathcal{C}}_{A}$ and is therefore relatively less attracting when the switching dynamics begin.

Fig.~\ref{fig:ResidenceTimes_xcen_50_65_rho_042_02}~(a) illustrates that for Case 1 there are a number of less strongly defined branches of exponentially distributed points. The two most well defined branches on the left hand side of this figure correspond to scenarios where the state of the closed-loop RC completes one or two loops (or partial loops) around $\mathcal{C}_{A}$, $\mathcal{C}_{B}$, or rapidly switches between $\tilde{\mathcal{C}}_{A}$ and $\tilde{\mathcal{C}}_{B}$. The well defined branch of orange points on the right hand side of this figure corresponds to the significantly longer amounts of time that the state of the closed-loop RC spends in $\tilde{\mathcal{C}}_{B}$ like in the example shown in Fig.~\ref{fig:xtraj_examples_CI_xcen_65_50_rho_020_042_}~(a) where the state of the closed-loop RC is in $\tilde{\mathcal{C}}_{B}$ from $t \approx 1295$ to $1390$. From further analysis we find that by increasing the number of bins, the cloud of points in the middle Fig.~\ref{fig:xtraj_examples_CI_xcen_65_50_rho_020_042_}~(a) corresponds to scenarios where the state of the closed-loop RC completes several loops (or partial loops) around $\mathcal{C}_{A}$ and $\mathcal{C}_{B}$. However, by increasing the number of bins we also find that this results in an increasingly large accumulation of points at the bottom of these branches which, in our opinion, diminishes the clarity of the message behind this figure and for that reason we do not present a version of Fig.~\ref{fig:xtraj_examples_CI_xcen_65_50_rho_020_042_}~(a) with a larger number of bins.

\subsubsection{Escape times}

\begin{figure}
    \centering
    \includegraphics[width=.95\textwidth]{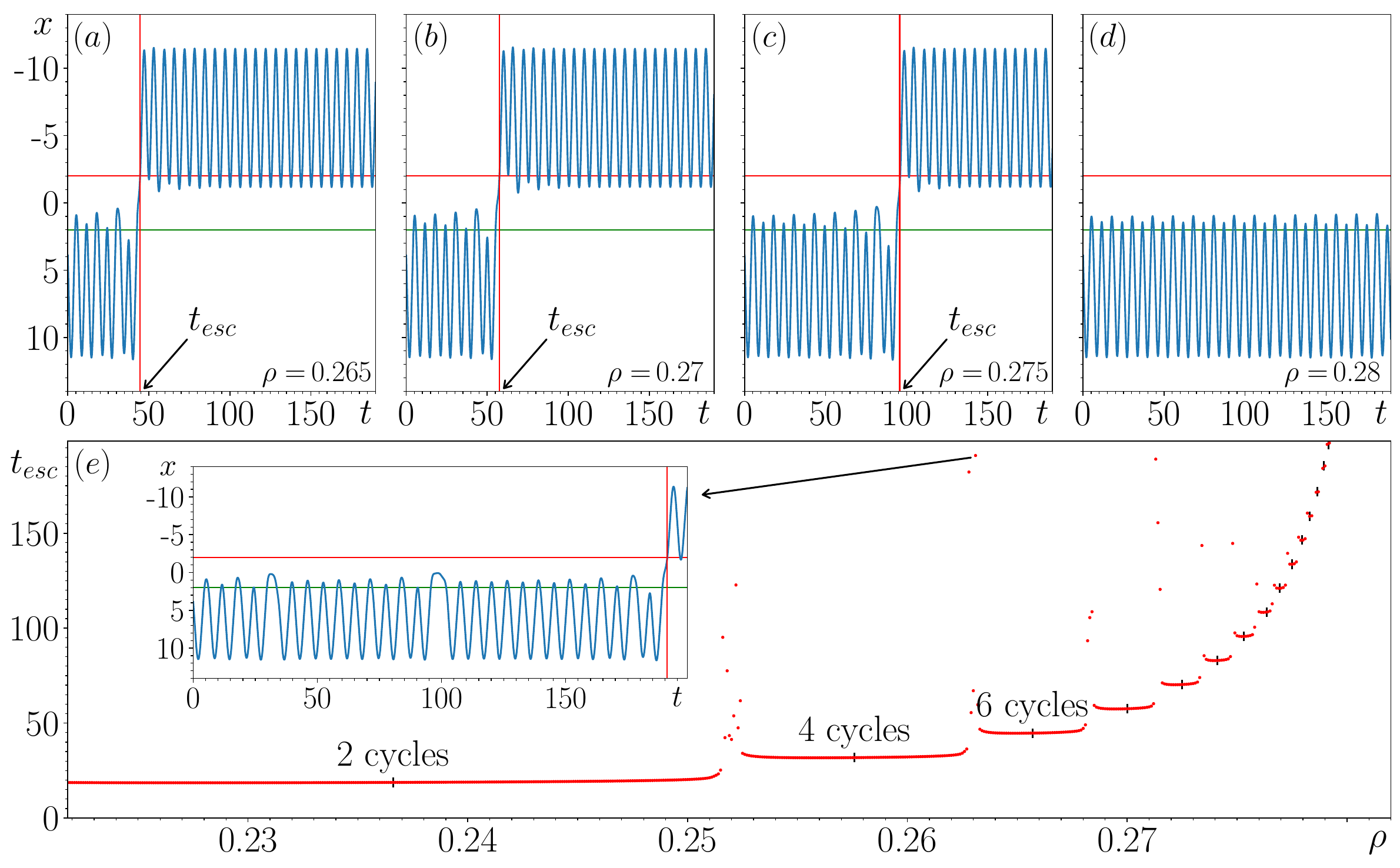}
    \caption{Illustrations of the closed-loop RC's dynamics when initialised from $\boldsymbol{r}_{(\mathcal{C}_{A})}(t_{train})$ in terms of $x(t)$ when trained for $x_{cen}=6.5$ and $\rho=0.265$ in panel (a), $\rho=0.27$ in panel (b), $\rho=0.275$ in panel (c), and $\rho=0.28$ in panel (d). In panel (e) we plot the values of $t_{esc}$ for $0.2218 \leq \rho \leq 0.28$, in the inset we plot how $x(t)$ behaves between one of the steps of the staircase-like structure seen in panel (e).}
    \label{fig:xtraj_examples_CI_xcen_65_rho_0265_027_0275_028_}
\end{figure}

The purpose of this section is to provide further insight to the interesting transient dynamics associated with $\hat{\mathcal{C}}_{A}$ becoming unstable when $x_{cen}=6.5$ as discussed in Sec.~\ref{sec:res_xcen65}. Using the transition detection algorithm we calculate the time it takes for the closed-loop RC to escape from transient behaviour when initialised from $\boldsymbol{r}_{(\mathcal{C}_{A})}(t_{train})$, we denote this duration of time as $t_{esc}$. We investigate the relationship between $\rho$ and $t_{esc}$ for values of $\rho$ when no switching dynamics occur, for $0.2218 \leq \rho \leq 0.28$. In panels (a)-(d) of Fig.\,\ref{fig:xtraj_examples_CI_xcen_65_rho_0265_027_0275_028_} we plot time series of the reconstructed $x(t)$ variable at particular $\rho$ values. We use the same $\alpha$ and $\beta$ thresholds as in the previous section, indicated by the red and green horizontal lines. The vertical red line depicts the detected value of $t_{esc}$. 

Figs.\,\ref{fig:xtraj_examples_CI_xcen_65_rho_0265_027_0275_028_}~(a)-(d) indicate that as $\rho$ increases and approaches the point at which $\hat{\mathcal{C}}_{A}$ becomes stable, then $t_{esc}$ naturally increases. However, panel (e) shows that while $t_{esc}$ increases as $\rho$ increases, $t_{esc}$ increases in a non-trivial staircase-like manner where the length of each successive step decreases as $\rho$ increases. 

For instance, $t_{esc}$ is shown here to be relatively small and increasing at a relatively slow rate for $\rho < 0.251$, however by increasing $\rho$ to $0.253$ this results in nearly a two-fold increase in $t_{esc}$, but for $0.251 < \rho < 0.253$ we find relatively large values and large variations in the values of $t_{esc}$ where $\hat{\mathcal{C}}_{A}$ appears to almost regain stability. The inset plot shows one of these relatively long transients between the steps at $\rho 0.2631$. Here we see from the change in time where local minima occur, the state of the closed-loop RC almost escapes from this transient activity at $t \approx 32$, as it does so for smaller $\rho$ values, and again at $t \approx 100$. This change in time is indicative of the state of the closed-loop RC approaching the saddle point on the unstable $\hat{\mathcal{C}}_{A}$ but fails to cross its separatrix. It is only for $t \approx 191$ that the RC escapes from the transient.

As indicated in panel (e), at each successive step along this staircase, the state of the closed-loop RC completes two additional cycles about the unstable $\hat{\mathcal{C}}_{A}$ before escaping. While the calculated values of $t_{esc}$ depend on the choice of initial condition relative to the point of escape on the unstable $\hat{\mathcal{C}}_{A}$, this behaviour of completing two additional cycles at each successive step may be more strongly dependent on the nature of $\hat{\mathcal{C}}_{A}$ prior to becoming unstable as it exists as a period-2 limit-cycle (indicated by the two global maxima Figs.\,\ref{fig:xtraj_examples_CI_xcen_65_rho_0265_027_0275_028_}~(d), also seen in Fig.~\ref{fig:bifplot_xcen6.5_}~(e) albeit barely visible). Our results suggest that by increasing $\rho$ the saddle point on the unstable $\hat{\mathcal{C}}_{A}$ moves in a way that the state of the closed-loop RC needs to complete an additional round trip about the unstable period-2 nature of $\hat{\mathcal{C}}_{A}$ until it reaches the point of escape.

\section{Discussion}

In this paper we explore how switching dynamics emerge in a dynamical system in the form of a RC when trained to achieve multifunctionality by solving the seeing double problem. This problem involves training the open-loop RC in Eq.~\eqref{eq:ListenRes} to reconstruct a coexistence of two circular orbits $\mathcal{C}_{A}$ and $\mathcal{C}_{B}$. We find that as $\mathcal{C}_{A}$ and $\mathcal{C}_{B}$ are moved closer together, the state of the closed-loop RC (Eq.~\eqref{eq:PredRes}) begins to switch between what appear to be metastable states that resemble trajectories around regions of $\mathbb{P}$ associated with $\mathcal{C}_{A}$ and $\mathcal{C}_{B}$. To be more specific, we find that these switching dynamics occur just before $\mathcal{C}_{A}$ and $\mathcal{C}_{B}$ touch (as shown in Fig.~\ref{fig:bifplot_xcen6.5_} for $x_{cen}=6.5$), as they touch (as shown in Fig.~\ref{fig:bifplot_xcen5.0_} for $x_{cen}=5.0$), and after they touch (as shown in Fig.~\ref{fig:bifplot_xcen2.0_} for $x_{cen}=2.0$) whereby there is an overlap between $\mathcal{C}_{A}$ and $\mathcal{C}_{B}$. However, as shown in Fig.~\ref{fig:bifplot_xcen3.5_}, there is an intermediary regime whereby after $\mathcal{C}_{A}$ and $\mathcal{C}_{B}$ touch and begin to overlap (for $x_{cen}=3.5$) the RC recovers its ability to achieve multifunctionality and does not succumb to these switching dynamics. It is only after there is a sufficiently large amount of overlap between $\mathcal{C}_{A}$ and $\mathcal{C}_{B}$ (for $x_{cen}=2.0$) that the switching dynamics reappear.

Our results also shed further light on the key role played by $\rho$ in this RC design and its connection to the concept of memory in terms of how the larger the value of $\rho$ the greater the influence of previous states on the current state of the RC. What our results indicate is that if the orbits are close to touching each other like for $x_{cen}=6.5$ or touch each other at only one point when $x_{cen}=5.0$, this requires the RC to place a greater weight on previous states (i.e. large $\rho$) in order to achieve multifunctionality as the dynamics nearby these touching regions are quite similar. On the other hand if the orbits overlap and touch each other in two locations that are sufficiently far but not too far apart, like for $x_{cen}=3.5$, then the RC does not need to place such a large weight on previous states in order to achieve multifunctionality. However once there is a larger amount of overlap between the orbits, like for $x_{cen}=2.0$, then the RC needs to place greater weight on previous states in order to achieve multifunctionality once again.

It is also worth noting that in panel (e) of Figs.\,\ref{fig:bifplot_xcen6.5_}-\ref{fig:bifplot_xcen2.0_}, prior to $\hat{\mathcal{C}}_{A}$ or $\hat{\mathcal{C}}_{B}$ becoming unstable as $\rho$ decreases there is a noticeable difference the obtained values for $x_{m}$ and the corresponding true values. This is most evident in panels (a)-(d) of Fig.~\ref{fig:bifplot_xcen3.5_} where we see $\hat{\mathcal{C}}_{A}$ and $\hat{\mathcal{C}}_{B}$ stretched towards larger positive and negative values of $x$ respectively. As $x_{cen}$ is decreased further this effect appears to becomes more and more noticeable. A similar sequence of events was shown to occur in Figs.~14, 15, and 21 in \cite{flynn2023seeingdouble} where, for $x_{cen}=0$, as $\rho$ decreases $\hat{\mathcal{C}}_{A}$ and $\hat{\mathcal{C}}_{B}$ are deformed in a similar way. This particular deformation may occur due to the design of $\textbf{W}_{in}$ as each neuron receives input from only one component of the driving input signal because each row contains only one nonzero element, therefore as $\rho$ decreases this increases the influence of the input and this may increase the likelihood that the resulting dynamics of the closed-loop RC are stretched along the $y=x$ and $y=-x$ diagonals. However, in order to provide a more rigorous answer this requires conducting an extensive analysis across several random realisations of $\textbf{M}$ and $\textbf{W}_{in}$ and testing whether such a deformation effect persists when using different design principles to construct $\textbf{M}$ and $\textbf{W}_{in}$. We believe that such an investigation is highly worthwhile to conduct and is better suited to appear in a paper where this is the main focus.

From closer inspection of the transitions between these metastable states, that we refer to as $\tilde{\mathcal{C}}_{A}$ and $\tilde{\mathcal{C}}_{B}$, we find that there is a common sequence of events that occurs in each case in order to produce the switchings between $\tilde{\mathcal{C}}_{A}$ and $\tilde{\mathcal{C}}_{B}$. Starting from a set of training parameters where the closed-loop RC achieves multifunctionality, we track how the dynamics of $\hat{\mathcal{C}}_{A}$ and $\hat{\mathcal{C}}_{B}$ change with respect to changes in $\rho$, the spectral radius of the RC's internal connectivity matrix. We find that by decreasing $\rho$ from the point where $\hat{\mathcal{C}}_{A}$ and $\hat{\mathcal{C}}_{B}$ coexist and resemble $\mathcal{C}_{A}$ and $\mathcal{C}_{B}$, there is a value of $\rho$ where, for instance, $\hat{\mathcal{C}}_{A}$ collides with a nearby saddle and becomes unstable but there still exists some transient dynamics that the state of the closed-loop follows when initialised from a point on the previously stable $\hat{\mathcal{C}}_{A}$. Then, by further decreasing $\rho$ we find that there is a value of $\rho$ where $\hat{\mathcal{C}}_{B}$ also becomes unstable by colliding with a nearby saddle. However when $\hat{\mathcal{C}}_{B}$ becomes unstable there is a new attractor born that facilitates the switching dynamics between the metastable states, $\tilde{\mathcal{C}}_{A}$ and $\tilde{\mathcal{C}}_{B}$, mentioned earlier. To be more specific, a trajectory on this new attractor consists of two regions of convergent flow where the trajectory inside these regions resembles a trajectory around $\mathcal{C}_{A}$ and $\mathcal{C}_{B}$, and a divergent flow whereby the state of the closed-loop RC switches from one region of convergent flow to the other. 

We also investigate the long-term behaviour of some of these new attractors that are born during the sequence of events discussed above. We integrated the closed-loop RC forwards in time until we obtain 10,000 transitions between $\tilde{\mathcal{C}}_{A}$ and $\tilde{\mathcal{C}}_{B}$ for the chaotic attractors illustrated in Fig.~\ref{fig:bifplot_xcen6.5_}~(c), denoted as Case 1, and in Fig.~\ref{fig:bifplot_xcen5.0_}~(d), denoted as Case 2. We construct an algorithm based on the concept of a non-ideal relay to determine the time of transition between $\tilde{\mathcal{C}}_{A}$ and $\tilde{\mathcal{C}}_{B}$. In Fig.~\ref{fig:xtraj_examples_CI_xcen_65_50_rho_020_042_} we provide an example of the transitions times detected by this algorithm. Interestingly, by computing the probability density of residence times in $\tilde{\mathcal{C}}_{A}$ and $\tilde{\mathcal{C}}_{B}$ we obtain several branches of exponentially distributed points as shown in Fig.~\ref{fig:ResidenceTimes_xcen_50_65_rho_042_02}. From closer inspection we find that each of these branches correspond to scenarios where the state of the closed-loop RC completes a given number of loops or partial loops around $\mathcal{C}_{A}$ and $\mathcal{C}_{B}$. 

We remark that while these switching dynamics are found for a particular random realisation of $\textbf{M}$ and $\textbf{W}_{in}$ (the internal and input connectivity matrices), the results presented in this paper are not solely dependent on these particular weights as from further experiments not shown here we see similar behaviour emerging. Furthermore, there is a noticeable imbalance in the behaviour of $\hat{\mathcal{C}}_{A}$ and $\hat{\mathcal{C}}_{B}$ despite the symmetry present in the training data. We believe that this is due to the particular random realisation of the $\textbf{M}$ and $\textbf{W}_{in}$ matrices happening to favour the reconstruction of one orbit over the other at particular parameter settings. From further analysis (also not shown in this paper), we find some small differences in the values of $\rho$ and the order of when $\hat{\mathcal{C}}_{A}$ and $\hat{\mathcal{C}}_{B}$ become unstable for different realisations of $\textbf{M}$ and $\textbf{W}_{in}$. As a further point, while the switching dynamics are induced by moving $\mathcal{C}_{A}$ and $\mathcal{C}_{B}$ closer together, it is still possible for switching dynamics to emerge between a reconstructed attractor and untrained attractors (attractors that the closed-loop RC produces that was not present during the training), or between the attractor a RC with symmetry is trained to reconstruct and its mirrored counterpart as shown in Fig.~2 in \cite{herteux2020Symm}. We suspect that when there is a competition between attractors, be it attractors that are manually moved closer together or attractors that compete with their mirrored counterpart or other untrained attractors, this sequence of attractors becoming unstable combined with the constraint that the RC is prohibited from exhibiting globally unstable dynamics (due to the choice of activation function) in turn creates a new attractor that is composed of different metastable states that in turn produces these switching dynamics.

Out of the many examples of routes to metastable dynamics discussed in \cite{rossi2024metastable}, there are a number of similarities between the results presented in this paper and phenomena such as chaotic itinerancy and heteroclinic cycles. In the case of chaotic itinerancy, which describes a switching process whereby the state of an autonomous dynamical system switches between several `attractor ruins' or `quasi-attractors' (these were previously coexisting attractors that retain much of their original features except trajectories on these quasi-attractors leak into each other) in our case these quasi-attractors are described as the metastable states $\tilde{\mathcal{C}}_{A}$ and $\tilde{\mathcal{C}}_{B}$. In terms of heteroclinic cycles, this typically occurs when the unstable manifold of one saddle coincides with a stable manifold of the other saddle which in our case these saddles would be the chaotic transients associated with $\mathcal{C}_{A}$ and $\mathcal{C}_{B}$. However further work is required in order to determine which of these phenomena our results are most closely aligned with. Furthermore, a similar route to chaotic behaviour has been observed in the past by \cite{grebogi1985super} whereby when two unstable orbits move towards each other by changing a parameter in the system, they coalesce at a bifurcation point and subsequently disappear, however after the bifurcation a chaotic transient is produced which persists for parameter values far beyond the bifurcation point. In our case we have one stable attractor and an unstable orbit/relatively long transient in the closed-loop RC that as $\rho$ is varied there is a bifurcation where the stable attractor becomes unstable and a new attractor is born which, depending on the circumstances, is either a chaotic attractor or limit cycle. Moreover, there is a valid reason why there is no transient produced after the second attractor becomes unstable. Due to the design of this closed-loop RC it is prevented from ever becoming globally unstable and since there is no other stable attractor present in the closed-loop RC when the second attractor becomes unstable there is no option but for there to be a stable attractor born through these sequence of attractors becoming unstable.

While the routes to metastable behaviour mentioned above are well-studied phenomena they only arise in certain circumstances. Rather than relying on there being a parameter in a dynamical system that so happens to produce these switching dynamics, the major advantage of the multifunctional reservoir computing setup studied in this paper is that we are able to systematically induce these switching dynamics by adjusting the location of $\mathcal{C}_{A}$ and $\mathcal{C}_{B}$. Moreover, by using this framework of multifunctionality its possible to study how these switching dynamics emerge in more complicated scenarios where the RC is trained to, for instance, reconstruct a coexistence of tori, or chaotic attractors, or a combination of different types of attractors. Given the rich variety of interesting dynamics that we see arise when training the RC to reconstruct a coexistence of two circular orbits we expect that in these more complicated scenarios there are even more interesting dynamics waiting to be explored, we leave such an investigation for future work.




As a final comment, the work presented throughout this paper highlights the importance of studying the behaviour of saddles and the bifurcations which take place as a RC, or any dynamical systems based machine learning approach, is trained to solve a given task. As strongly emphasised in \cite{Sussillo13_PCA}, in order to open the black-box of machine learning approaches, it is necessary that we improve our understanding of the interaction between stable and unstable dynamics and pay closer attention to the influence of saddles that are present in the system. 





\section*{Conflict of Interest Statement}

The authors declare that the research was conducted in the absence of any commercial or financial relationships that could be construed as a potential conflict of interest.

\section*{Author Contributions}

AF: Writing – review \& editing, Conceptualization, Methodology; AA: Writing – review \& editing.


\section*{Acknowledgments}
We would like to thank Aravind Kumar and our Applied Mathematics colleagues at the School of Mathematical Sciences, in particular Andrew Keane, Pierce Ryan, Serhiy Yanchuk, and Sebastian Wieczorek, for their influential conversations and input when discussing the contents of this paper. This work was partly supported by the Deutsche Forschungsgemeinschaft Project No. 411803875 and PIK Werkvertrag 2023-0336.


\section*{Data Availability Statement}
The datasets generated for this study are available from the corresponding author upon reasonable request.

\bibliographystyle{unsrtnat}






\end{document}